\documentclass{amsart}
\usepackage{amsxtra}
\theoremstyle{plain}

\newcommand{\cleqn}{\setcounter{equation}{0}}
\newcommand{\clth}{\setcounter{theorem}{0}}
\newcommand {\sectionnew}[1]{\section{#1}\cleqn\clth}
\newcommand{\nn}{\hfill\nonumber}
\newtheorem{theorem}{Theorem}[section]
\newtheorem{lemma}[theorem]{Lemma}
\newtheorem{definition-theorem}[theorem]{Definition-Theorem}
\newtheorem{proposition}[theorem]{Proposition}
\newtheorem{corollary}[theorem]{Corollary}
\newtheorem{definition}[theorem]{Definition}
\newtheorem{example}[theorem]{Example}
\newtheorem{remark}[theorem]{Remark}
\newtheorem{notation}[theorem]{Notation}
\newcommand \bth[1] { \begin{theorem}\label{t#1} }
\newcommand \ble[1] { \begin{lemma}\label{l#1} }

\newcommand \bpr[1] { \begin{proposition}\label{p#1} }
\newcommand \bco[1] { \begin{corollary}\label{c#1} }
\newcommand \bde[1] { \begin{definition}\label{d#1}\rm }
\newcommand \bex[1] { \begin{example}\label{e#1}\rm }
\newcommand \bre[1] { \begin{remark}\label{r#1}\rm }

\newcommand \bnota[1] { \begin{notation}\label{n#1}\rm }
\newcommand {\eth} { \end{theorem} }
\newcommand {\ele} { \end{lemma} }

\newcommand {\epr} { \end{proposition} }
\newcommand {\eco} { \end{corollary} }
\newcommand {\ede} { \end{definition} }
\newcommand {\eex} { \end{example} }
\newcommand {\ere} { \end{remark} }

\newcommand {\enota} { \end{notation} }
\newcommand \thref[1]{Theorem \ref{t#1}}
\newcommand \leref[1]{Lemma \ref{l#1}}
\newcommand \prref[1]{Proposition \ref{p#1}}
\newcommand \coref[1]{Corollary \ref{c#1}}
\newcommand \deref[1]{Definition \ref{d#1}}
\newcommand \exref[1]{Example \ref{e#1}}
\newcommand \reref[1]{Remark \ref{r#1}}
\newcommand \lb[1]{\label{#1}}

\def \A  {{\mathcal{A}}}           

\def \C  {{\mathcal{C}}}

\def \Ocal  {{\mathcal{O}}}

\def \Zcal  {{\mathcal{Z}}}
\def \Xcal  {{\mathcal{X}}}

\def \Wcal {{\mathcal W}}

\def \Jcal  {{\mathcal{J}}}

\def \al {\alpha}

\def \Ga {\Gamma}

\def \sig {\sigma}

\def \sig{\sigma}

\def \ra  {\rightarrow}           
\def \lra {\longrightarrow}

\def \ol {\overline}

\def \hs {\hspace{.2in}}

\def \id { {\mathrm{id}} }
\def \Stab { {\mathrm{Stab}} }

\def \Kvv  { {K_{(v_1, v_2)}} }
\def \Lvv  { {L_{(v_1, v_2)}} }

\def \Qvv { {Q_{(v_1, v_2)}} }
\def \Rvv {{R_{(v_1, v_2)}}}

\def \Aavv { A_1(v_1, v_2) }
\def \Abvv { A_2(v_1, v_2) }

\def \Cbvv { C_2(v_1, v_2) }

\def \Ad { {\mathrm{Ad}} }

\def \k {\mathrm{k}}

\DeclareMathOperator \Aut { {\mathrm{Aut}} }

\def \dw {\dot{w}}
\def \dv {\dot{v}}
\def \du {\dot{u}}

\def \dx {\dot{x}}

\def \vvprime {[v_1, v_2]_{\A, \C}}
\def \vvAC {[v_1, v_2]_{\A, \C}}
\def \vvACm {[v_1, v_2]_{\A, \C}^{-}}
\def \wg {w_{0, \Gamma_1}}
\def \wc {w_{0, C_1}}
\def \dwg {\dot{w}_{0, \Gamma_1}}
\def \dwc {\dot{w}_{0, C_1}}

\def \wcc {w_{0, C}}

\begin{document}
\title[Partitions of the wonderful compactification]
{Partitions of the wonderful group compactification}
\author[Jiang-Hua Lu]{Jiang-Hua Lu}
\address{
Department of Mathematics   \\
The University of Hong Kong \\
Pokfulam Road               \\
Hong Kong}
\email{jhlu@maths.hku.hk}
\author[Milen Yakimov]{Milen Yakimov}
\address{
Department of Mathematics \\
University of California \\
Santa Barbara, CA 93106, U.S.A.}
\email{yakimov@math.ucsb.edu}
\date{}
\begin{abstract}
We define and study a family of partitions of the wonderful 
compactification $\ol{G}$ of a semi-simple algebraic group 
$G$ of adjoint type.
The partitions are obtained from subgroups of $G \times G$
associated to
triples $(A_1, A_2, a)$, where $A_1$ and $A_2$ are subgraphs 
of the Dynkin graph $\Ga$ of $G$ and $a \colon A_1 \ra A_2$
is an isomorphism. The partitions of $\ol{G}$ of Springer and Lusztig
correspond respectively to the triples $(\emptyset, \emptyset, \id)$ and 
$(\Ga, \Ga, \id).$
\end{abstract}
\maketitle
\sectionnew{Introduction}\lb{intro}
Let $G$ be a connected semi-simple algebraic group over an algebraically 
closed
field $k$. 
De Concini and Procesi \cite{DP, DS} constructed a wonderful compactification
$\ol{G}$ of $G$, which is a smooth irreducible $(G \times G)$-variety 
with finitely many
$(G \times G)$-orbits. 
Let $G_{{\rm diag}}$ be the diagonal subgroup of $G \times G$.
In his study of parabolic character sheaves on $\ol{G}$ in \cite{L1, L2}, 
Lusztig
introduced (by an inductive procedure) a partition of $\ol{G}$ by 
finitely many 
$G_{{\rm diag}}$-stable pieces. The closure of a $G_{{\rm diag}}$-stable
piece was shown by X. He \cite{H2} to be a union of such pieces. 
Let $B$ be a Borel subgroup of $G$. Then $\ol{G}$ is also partitioned 
into finitely many
$(B \times B)$-orbits. The $(B \times B)$-orbits in $\ol{G}$,
as well as their closures, were studied by T. Springer in \cite{Sp}. 
In \cite{H2},
X. He gave a second description of Lusztig's $G_{{\rm diag}}$-stable 
pieces using
$(B \times B)$-orbits in $\ol{G}$, which then enabled him to give \cite{H3}
an equivalent definition of Lusztig's character sheaves on $\ol{G}$.
Further properties and applications of the $G_{{\rm diag}}$-stable pieces 
were obtained
by X. He and J. F. Thomsen in \cite{H1, H3, HT1} and T. A. Springer in \cite{Sp1}.

Both $G_{{\rm diag}}$ and $B \times B$ are special examples of subgroups 
$R_\A$ of
$G \times G$ associated to triples $(A_1, A_2, a)$, where
$A_1$ and $A_2$ are subgraphs of the Dynkin graph $\Ga$ of $G$, and $a$ 
is an isomorphism
from $A_1$ to $A_2$. If $P_{A_1}$ and $P_{A_2}$ are the standard 
parabolic subgroups of $G$ 
corresponding to $A_1$ and $A_2$ respectively, then, roughly speaking, 
$R_\A$ is a subgroup of $P_{A_1} \times P_{A_2}$, obtained by 
identifying the Levi subgroups 
of $P_{A_1}$ and $P_{A_2}$ via the map $a$.
The precise definition of $R_\A$ is given in $\S$\ref{quad}. For example, 
every stabilizer 
subgroup of $G \times G$ in $\ol{G}$ is conjugate to a group of this form. 
Moreover, $G_{{\rm diag}}$ is associated to the triple $(\Ga, \Ga, {\rm 
id})$ and $B \times B$
to the triple 
$(\emptyset, \emptyset, {\rm id})$, where $\emptyset$ is the empty set. 

In this paper, 
for any subgroup $R_\A$ of $G \times G$ associated to a triple $(A_1, A_2, a)$,
we study a partition of $\ol{G}$ into finitely many $R_\A$-stable pieces indexed by a subset of
the Weyl group of $G \times G$. 
Our definition of the $R_\A$-stable pieces is based on 
our earlier paper \cite{LY} on
$(R_\A, R_\C)$-double cosets in $G \times G$ for subgroups $R_\A$ and $R_\C$
associated to any pair of triples $(A_1, A_2, a)$ and $(C_1, C_2, c)$.
We give two additional descriptions of the $R_\A$-stable pieces, 
which for the case of $R_\A = G_{\rm diag}$, reduce to Lusztig's 
inductive description in \cite{L1, L2}
and He's description in \cite{H2} using $(B \times B)$-orbits. 
In particular, we show that the $R_\A$-stable pieces are smooth,
irreducible, locally closed
subsets of $\ol{G}$, fibering over flag varieties of 
Levi subgroups of $G$. We also show that the closure in $\ol{G}$ 
of an $R_\A$-stable piece is a union of such pieces.
We describe the combinatorics for the closures of the $R_\A$-stable 
pieces, generalizing both the result of He \cite{H2} for the $G_{{\rm diag}}$-stable pieces
and that of Springer for the $(B \times B)$-orbits. The closure
relations of the $R_\A$-stable pieces are expressed in terms of intersections of the closures
with the unique closed $(G \times G)$-orbit in $\ol{G}$. 

Our
motivation for studying the $R_\A$-stable pieces in $\ol{G}$ for arbitrary triples 
$(A_1, A_2, a)$
comes from Poisson geometry. 
In \cite{LY3}, we study a class of
Poisson structures on $\overline{G}$. Special examples of these
Poisson structures are
induced by the Belavin--Drinfeld
$r$-matrices \cite{BD}. The triples $(A_1, A_2, a)$ needed there
are precisely the Belavin--Drinfeld triples for the $r$-matrices.
The $R_\A$-stable pieces in $\ol{G}$ as well as their closures are Poisson
subvarieties of $\ol{G}$ for the corresponding Poisson structures. To understand
these Poisson structures, one needs to first understand the geometry of the
$R_\A$-stable pieces.

In $\S$\ref{setup} - $\S$\ref{closures-GG}, we in fact assume that 
$G_1$ and $G_2$ are any two
reductive algebraic groups over an algebraically closed field, and that
$R_\A$ and $R_\C$ are subgroups of $G_1 \times G_2$ associated to 
two triples
$(A_1, A_2, a)$ and $(C_1, C_2, c)$
for $G_1 \times G_2$. The precise definitions of
the subgroups $R_\A$ and $R_\C$
are given in $\S$\ref{quad}. Each pair $(R_\A, R_\C)$ of such subgroups
gives rise to a decomposition of $G_1 \times G_2$ into $(R_\A, R_\C)$-stable
subsets of the form 
$[v_1, v_2]_{\A, \C}$, where $(v_1, v_2)$ runs over a subset of the Weyl group
for $G_1 \times G_2$ (see \eqref{partAC} for detail). 
In $\S$\ref{inductive-description}, we give a description of
the set  $[v_1, v_2]_{\A, \C}$ as iterated fiber bundles.
The closures of the sets $[v_1, v_2]_{\A, \C}$ in $G_1 \times G_2$ are described in 
$\S$\ref{closures-GG}. 
We point out that in a recent paper
\cite{H4}, 
X. He studies the $(R_\A, R_\C)$-stable pieces in $G_1 \times G_2$ 
from the point of view of Coxeter groups. 
In $\S$\ref{minus-strata} and $\S$\ref{proof_th}, the 
results in $\S$\ref{closures-GG} are used to prove our main theorems on the 
$R_\A$-stable pieces in $\ol{G}$ for a semi-simple algebraic group $G$ of adjoint type. 
Precise statements of our results are summarized in $\S$\ref{setup}. In 
the appendix we collect a few facts on the Bruhat order 
on Weyl groups that we use in the paper.

{\bf Acknowledgements.} 
We would like to thank Sam Evens and Xuhua He for helpful answers to
our questions. We also thank the referees whose comments
helped us to improve the exposition.
\sectionnew{Notation and statements of results}\label{setup}
\subsection{Admissible quadruples} 
\label{quad}
For $i = 1, 2$, let $G_i$ be a connected reductive algebraic group 
over an algebraically closed field $\k$. Let $B_i$ and $B_i^-$ be 
a fixed pair of opposite Borel subgroups of $G_i$ with $U_i$ being their
respective uniradicals. 
Set $T_i = B_i \cap B_i^-$, 
and let
$\Gamma_i$ be the set of simple roots determined by $(B_i, T_i)$.
For 
$\al \in \Ga_{i}$, denote by $U^\al_i$ the one-parameter 
unipotent subgroup of $G_i$ defined by $\al.$ 
For a subset $A_i$ of $\Gamma_i$, let $P_{A_i}$ and $P_{A_i}^-$ 
be the standard parabolic subgroups of $G_i$ containing respectively 
$B_i$ and $B_i^-$. Let $M_{A_i}=P_{A_i} \cap P_{A_i}^-$
be the common Levi factor of $P_{A_i}^\pm$, and let 
$Z_{A_i}$ be the center of $M_{A_i}$.
The unipotent radicals of $P_{A_i}$ and $P_{A_i}^-$ will be denoted by
$U_{A_i}$ and $U_{A_i}^-$ respectively. 
Let $W_i$ be the Weyl group of $\Ga_i$
and $W_{A_i}$ the subgroup
of $W_i$ generated by reflections
defined by simple roots in $A_i$.
Let $W_i^{A_i}$ and ${}^{A_i} \! W_i$ be the sets of minimal length 
representatives of cosets from $W_i/W_{A_i}$ and 
$W_{A_i} \backslash W_i$ respectively. For each $w_i \in W_i$, we also fix
a choice $\dw_i$ of a representative of $w_i $ in the normalizer of $T_i$
in $G_i$. The length function on $W_i$ will be denoted by $l$ and the Bruhat order by
$\leq$.
If a group $G$ acts on a set $X$,  
$g_\cdot x$ denotes the action of $g \in G$ on $x \in X$. 
For an element $g \in G$, the map
$G \to G: h \mapsto g h g^{-1}$ will be denoted by $\Ad_g$.
The identity element of a group will be denoted by $e$ or $1$.

For subsets $A_i$ of the Dynkin graphs $\Gamma_i$, $i=1,2$,
we call a bijective map $a \colon A_1 \to A_2$ an isomorphism, 
if it preserves the type of each arrow. 

\bde{admis} An admissible 
quadruple  for 
$G_1 \times G_2$ is a quadruple  $\A=(A_1, A_2, a, K)$
consisting of subsets $A_1$ of $\Ga_1$ and $A_2$ of $\Ga_2$, an
isomorphism
$a \colon A_1 \to A_2$, and a closed subgroup $K$ of
$M_{A_1} \times M_{A_2}$ of the form
\begin{equation}
\label{K}
K = \{(m_1, m_2) \in M_{A_1} \times M_{A_2} \mid \theta_a (m_1Z_1) =
m_2Z_2\},
\end{equation}
where, for $i = 1, 2$,  $Z_i$ is a closed subgroup of $Z_{A_i}$ and 
$\theta_a: M_{A_1}/Z_1
\to M_{A_2}/Z_2$ is an isomorphism mapping $T_1/Z_1$ to $T_2/Z_2$ and 
$U_{1}^{\alpha}$ to $U_{2}^{a(\alpha)}$
for each $\alpha \in A_1$. Here we identify $U_{1}^{\alpha}$
and $U_{2}^{a(\alpha)}$ with their images in $M_{A_1}/Z_1$ and 
$M_{A_2}/Z_2$ respectively.
Given an admissible quadruple $\A = (A_1, A_2, a, K)$ of $G_1 \times G_2$,
define  
\begin{equation}
\label{Ra}
R_\A =
K (U_{A_1} \times U_{A_2}) \subset P_{A_1} \times P_{A_2}, \hspace{.3in}
R_\A^- = 
K (U_{A_1}^{-} \times U_{A_2}) \subset P_{A_1}^- \times P_{A_2}.
\end{equation}
\ede

Note that when $G_1 = G_2 = G$, the diagonal subgroup $G_{{\rm diag}}$ 
and $B \times B$ for a Borel subgroup $B$ are examples of the groups $R_\A$.
If further $G$ is of adjoint type, all
stabilizer subgroups of $G \times G$ in the De Concini--Procesi
\cite{DP} compactification $\ol{G}$ of $G$ are conjugate to 
groups of the type $R_\A^-$ (see \S \ref{wond}).

\subsection{An $(R_\A, R_\C)$-stable partition of $G_1 \times G_2$}
\label{G12part}
In \cite{LY} we obtained a 
classification  of $(R_\A, R_\C)$-double cosets 
of $G_1 \times G_2$ for two arbitrary admissible quadruples
$\A=(A_1, A_2, a, K)$ and $\C=(C_1, C_2, c, L)$ for $G_1 \times G_2$.
Given $v_1 \in W_1^{C_1}, v_2 \in {}^{A_2} \! W_2$, 
set
\begin{equation}
\label{C2stab}
C_2(v_1, v_2) =\{ \beta \in C_2 \mid
(v_{2}^{-1} a v_1 c^{-1})^n \beta \; 
\mbox{is defined and is in} \; C_2 \;
\mbox{for} \; n=1, 2, \ldots\}.
\end{equation}
In other words, $C_2(v_1, v_2)$ is the largest subset of $C_2$ 
that is stable under $v_{2}^{-1}a v_1 c^{-1}$. 
We proved \cite{LY}
that each $(R_\A, R_\C)$-double coset of $G_1 \times G_2$ is
of the form $R_\A(\dot{v}_1, \dot{v}_2 m_2) R_\C$ for some 
$v_1 \in W_1^{C_1}, v_2 \in {}^{A_2} \! W_2$ 
and $m_2 \in M_{\Cbvv}.$ Moreover, 
two such double cosets 
$R_\A(\dot{v}_1, \dot{v}_2m_2) R_\C$ and 
$R_\A( \dot{v}'_1, \dot{v}'_2m_2^\prime) R_\C$ coincide 
if and only if $v'_i=v_i$ for $i = 1, 2,$ and 
$m_2$ and $ m^\prime_2$ 
are in the same $(v_{2}^{-1} a v_1 c^{-1})$-twisted 
conjugacy class in $M_{\Cbvv},$ see \thref{LYmain1}
for details.

For $v_1 \in W_1^{C_1}$ and $v_2 \in {}^{A_2} \! W_2$, let 
\begin{equation}
\label{partAC}
[v_1, v_2]_{\A, \C} = R_\A ( v_1, v_2 M_{\Cbvv}) R_\C
\subset G_1 \times G_2.
\end{equation}
Then by the above result from \cite{LY}, 
we have the decomposition
\begin{equation}
\label{sqcupAC}
G_1 \times G_2 = \bigsqcup_{v_1 \in W_1^{C_1},
v_2 \in {}^{A_2} \! W_2} [v_1, v_2]_{\A, \C}.
\end{equation}
Here and below $\bigsqcup$ stands for disjoint union.
Note that  \eqref{sqcupAC} is constructed in such a way that 
the $(R_\A, R_\C)$-double cosets of $G_1 \times G_2$
corresponding to the same discrete parameters
$v_1 \in W_1^{C_1}$ and 
$v_2 \in {}^{A_2} \! W_2$  but 
possibly different continuous parameters $m_2 \in M_{C_2(v_1, v_2)}$
are put together in a single stratum.  Alternatively
we have the decomposition
\begin{equation}
\label{sqcupAC2}
(G_1 \times G_2)/R_\C = \bigsqcup_{v_1 \in W_1^{C_1},
v_2 \in {}^{A_2} \! W_2} [v_1, v_2]_{\A, \C}/R_\C
\end{equation}
of $(G_1 \times G_2)/R_\C$ into $R_\A$-stable subsets.

The main objects of study in this paper
are the sets $[v_1, v_2]_{\A, \C}$ 
for $v_1 \in W_1^{C_1},
v_2 \in {}^{A_2} \! W_2$.
We describe their geometry, as well as their closure relations. The
results are then applied to the wonderful group compactifications.

The following theorem summarizes our
results for the decompositions \eqref{sqcupAC} and \eqref{sqcupAC2}, 
see \coref{properties},
\prref{2dc-plus-2}, \prref{pr-shift}, and \thref{th-vvAC-closures}.

\bth{main1} Given any two admissible quadruples $\A$ and $\C$ for 
$G_1 \times G_2$, the following hold for every
$ v_1 \in W_1^{C_1}$ and $v_2 \in {}^{A_2} \! W_2$.

(i) $[v_1, v_2]_{\A, \C}$ is locally closed, smooth, and 
irreducible.
Its projection $[v_1, v_2]_{\A, \C}/R_\C$ to 
$(G_1 \times G_2) /R_\C$ fibers over the
flag variety $M_{A_1}/(M_{A_1} \cap P_{\Aavv})$
with fibers isomorphic to the product of $M_{\Cbvv}/Y_2$ 
and the affine space of 
dimension 
\[
\dim U_{A_1(v_1, v_2)} - \dim (U_1 \cap \Ad_{\dv_1}U_{C_1}) + l(v_2),
\]
where
$A_1(v_1, v_2) = v_1 c^{-1} C_2(v_1, v_2) \subset A_1,$ and
$Y_2 =\{m \in M_{C_2} \mid (e, m) \in L\} \subset Z_{C_2}$. 

(ii) Alternatively the set $[v_1, v_2]_{\A, \C}$ is given by
\begin{equation}
[v_1, v_2]_{\A, \C} = R_\A (B_1\times B_2)(v_1, v_2)R_\C.
\label{second}
\end{equation}

(iii) The Zariski closure of $[v_1, v_2]_{\A, \C}$ in $G_1 \times G_2$
consists of 
those $[w_1, w_2]_{\A, \C}$ with 
$w_1\in W_1^{C_1}$ and $w_2 \in {}^{A_2}\!W_2$ for which there
exist $x_1 \in W_{A_1}$ and $y_1 \in W_{C_1}$ 
such that 
\[
x_1 w_1 y_1 \leq v_1 \hspace{.2in} {\rm and} \hspace{.2in}
a(x_1) w_2 c(y_1) \leq v_2.
\]

(iv) If $\A^\prime = (A_1, A_2, a, K^\prime)$ and $\C^\prime = (C_1, C_2, c, L^\prime)$
are two other admissible quadruples containing the same triples $(A_1, A_2, a)$ and
$(C_1, C_2, c)$, then there exist $t_2, s_2 \in T_2$
such that
\[
[v_1, v_2]_{\A^\prime, \C^\prime} = (e, t_2) [v_1, v_2]_{\A, \C} (e, s_2),\hs
\forall v_1 \in W_1^{C_1}, v_2 \in {}^{A_2} \! W_2.
\]
\eth

\begin{example}\label{Bruhat} When
$\A=\C= (\emptyset, \emptyset, \id, T_1 \times T_2)$ so that
$R_\A = R_\C = B_1 \times B_2$, we have
\[
[v_1, v_2]_{\A, \C}= (B_1 \times B_2) (v_1, v_2) 
(B_1 \times B_2), \hspace{.2in} \forall (v_1, v_2) \in W_1 \times W_2.
\]
Thus \eqref{sqcupAC} reduces to the Bruhat decomposition
\[
G_1 \times G_2 =  \bigsqcup_{v_1 \in W_1, v_2 \in W_2} 
(B_1 \times B_2) (v_1, v_2) (B_1 \times B_2). 
\]
Part (iii) of \thref{main1} in this case is the 
well-known statement for the closures of Bruhat cells.
\end{example}
\subsection{Partitions of the wonderful group compactification} 
\label{wond}
Now we specialize to the case when $G_1$ and 
$G_2$ are both isomorphic to a connected semisimple 
algebraic group $G$ of 
adjoint type. All data for $G$ 
will be denoted as in \S \ref{quad}, omitting the index $i.$ 
In particular,
$B$ and $B^-$ will be two fixed opposite Borel subgroups of $G$, 
$T = B \cap B^-$,
and $\Gamma$ will 
be the set of simple roots determined by $(B, T)$.

For $J \subset \Ga$, let $\pi_J \colon M_J \to M_J/Z_J$ 
be the natural projection.
By abuse of notation, we will denote by $J$ 
the quadruple $(J, J, \id, L_J),$ where 
\[
L_J = \{(m_1, m_2) \in M_J \times 
M_J \mid \pi_J(m_1) = \pi_J(m_2)\};
\]
so
\begin{equation}
\label{RJ-minus}
R_J^- = L_J (U_J^- \times U_J) = 
\{(p_1, p_2) \in P_J^- \times P_J \mid 
\pi_J(p_1) = \pi_J(p_2)\}.
\end{equation}

Recall that the wonderful compactification $\ol{G}$ of $G$ is a smooth 
irreducible 
projective $(G \times G)$-variety containing $G$ as an open $(G \times G)$-orbit.
It was defined and studied 
(in the more general framework of symmetric spaces)
by De Concini and Procesi \cite{DP} in the complex case 
and by De Concini and Springer \cite{DS} for the general case of 
an algebraically closed field $\k.$
The $(G \times G)$-orbits on $\ol{G}$ are parameterized
by the subsets $J$ of $\Ga$. 
A base point $h_J$ of 
the $(G \times G)$-orbit corresponding to $J \subset \Gamma$
has stabilizer subgroup $R_J^-.$ 

The partitions \eqref{sqcupAC} induce partitions 
of $\ol{G}$ as follows. Given an 
arbitrary admissible quadruple $\A=(A_1, A_2, a, K)$ 
for $G \times G$, by setting
\begin{equation}
\label{wond_part}
[J, v_1, v_2]_\A = R_\A (B \times B) (v_1, v_2)_\cdot h_J, 
\quad J \subset \Gamma, 
v_1 \in W^J, v_2 \in {}^{A_2} \! W,
\end{equation}
we obtain the following partition of the wonderful compactification
(see $\S$\ref{vvACm-def} for detail)
\begin{equation}
\label{wond_part2}
\ol{G} = \bigsqcup_{J \subset \Gamma, 
v_1 \in W^J,
v_2 \in {}^{A_2} \! W} [J, v_1, v_2]_\A,
\end{equation}
cf. \eqref{sqcupAC2} and (ii) of \thref{main1}. We will refer to the sets
$[J, v_1, v_2]_\A$ in \eqref{wond_part2} as $R_\A$-stable pieces of $\ol{G}$.
If the subgroup $K$ 
in $\A$ is changed to $K^\prime$, the partition \eqref{wond_part2}
is changed by an overall left translation by $(e, t)$ for some $t \in T$
(see (iv) of \thref{main1} and \prref{pr-shift}). 

\bex{two-examples}
Let $\A$ be the trivial admissible quadruple 
$(\emptyset, \emptyset, {\rm Id}, T \times T)$.
Then $R_\A = B \times B$ and 
we recover Springer's partition \cite{Sp} of $\ol{G}$ by 
$(B\times B)$-orbits:
\begin{equation}
\label{Springer}
\ol{G} = \bigsqcup_{J \subset \Gamma, v_1 \in W^J, v_2 \in W} 
(B \times B)(v_1, v_2)_\cdot h_J.
\end{equation}
On the other hand,
let $\A$ be the quadruple 
$\Ga_{\rm diag}:=(\Gamma, \Gamma, \id, G_{\rm diag})$
where $G_{\rm diag}$ is the diagonal subgroup of $G \times G$.
By \cite{H2},  we recover Lusztig's 
partition \cite{L1, L2} of $\ol{G}:$
\begin{equation}
\label{Lusztig}
\ol{G} = \bigsqcup_{J \subset \Gamma, v_1 \in W^J} 
[J, v_1, 1]_{\Gamma_{\rm diag}}=G_{\rm diag} 
(B \times B)(v_1, 1)_\cdot h_J.
\end{equation}
\eex

For a general admissible quadruple $\A$ for $G \times G$, 
our partition \eqref{wond_part2} is a
discrete interpolation between Springer's partition \eqref{Springer}
and Lusztig's partition \eqref{Lusztig} of $\ol{G}$. 

The closure relations 
of the strata in \eqref{wond_part2} can be derived directly from part (iii) in \thref{main1}.
This is stated in \prref{withw0C}.
However, we give a more elegant description of the closures using the unique
closed $(G \times G)$-orbit in $\ol{G}$, namely the 
orbit $(G \times G)_\cdot h_\emptyset \cong G/B^- \times G/B$,
where $\emptyset$ is the empty subset of $\Ga$. 
More precisely, for 
any $J \subset \Ga$ and $v_1 \in W^J, \, v_2  \in{}^{A_2} \! W$,
let 
$\overline{[J, v_1, v_2]_\A}$ be the closure of 
$[J, v_1, v_2]_\A$ in $\overline{G}$, and set
\[
\partial_{\emptyset}\overline{[J, v_1, v_2]_\A}
= \overline{[J, v_1, v_2]_\A} \cap (G \times G)_\cdot h_{\emptyset}.
\]

The following is a summary of \thref{th-boundary} and \thref{th-third-descrip}.

\bth{bd-bd}
Let $\A$ be an admissible quadruple for 
$G \times G$. Let $I, J \subset \Gamma,$
$v_1 \in W^J, v_1^\prime \in W^I$, and 
$v_2, v_2^\prime \in  {}^{A_2} \! W$. Then

(i) $[I, v_1^\prime, v_2^\prime]_\A \subset \overline{[J, v_1, v_2]_\A}$ if and only if
$I \subset J$ and $[\emptyset, v_1^\prime, v_2^\prime]_\A \subset
\partial_\emptyset \overline{[J, v_1, v_2]_\A}$;

(ii) $[I, v_1^\prime, v_2^\prime]_\A \subset \overline{[J, v_1, v_2]_\A}$ if and only if
$I \subset J$ and 
 there exist
$x \in W_{A_1}$ and $z \in W_J$ such that $x v_1^{\prime} \geq v_1z$ and 
$a(x) v_2^\prime \leq v_2z.$
\eth
  
\bex{two-cases}
Consider the case when
$\A = (\emptyset, \emptyset, {\rm Id}, T \times T)$,  so that
$R_\A = B \times B$. \thref{bd-bd} implies that for any 
$I, J \subset \Ga$, $v_1 \in W^J, v_1^\prime \in W^I$,
and $v_2, v_2^\prime \in W$, 
\[
(B \times B)(v_1^\prime, v_2^\prime)_\cdot h_I \subset 
\overline{(B \times B)(v_1, v_2)_\cdot h_J}
\]
if and only if $I \subset J$ and 
there exists $z \in W_J$  such that $v_1^\prime \geq v_1 z$
and $v_2^\prime \leq v_2z$.
This description of the $(B \times B)$-orbit 
closures in $\overline{G}$, which is
simpler than what is given in \cite{Sp}, 
was independently obtained in
\cite{HT2}. 
When $\A = \Gamma_{\rm diag}$ as in \exref{two-examples}, 
\thref{bd-bd} implies that for any 
$I, J \subset \Ga$, $v_1 \in W^J, v_1^\prime \in W^I$,
\[
G_{{\rm diag}}(B \times B)(v_1^\prime, 1)_\cdot h_I \subset 
\overline{G_{{\rm diag}}(B \times B)(v_1, 1)_\cdot h_J}
\]
if and only if $I \subset J$ and 
there exists $x \leq z \in W_J$ such that $x v_1^\prime \geq v_1 z.$
It follows from \cite[Corollary 3.4]{H2} (see \leref{wwu} in the Appendix) that
the above description of the closures of the 
$G_{{\rm diag}}$-stable subsets
is the same as that given in \cite{H2} by X. He.
\eex
\sectionnew{$(R_\A, R_\C)$-double cosets in $G_1 \times G_2$
and their stablizer subgroups}
\label{stabAC}

In this section, we first recall some results from \cite{LY}.
  
\subsection{Classification of 
$(R_\A, R_\C)$-double cosets in $G_1 \times G_2$}
\label{double-cosets}
Fix two arbitrary admissible quadruples
$\A=(A_1, A_2, a, K)$ and $\C=(C_1, C_2, c, L)$ 
for $G_1 \times G_2.$
In \cite{LY} we obtained a classification of the
$(R_\A, R_\C)$-double cosets of $G_1 \times G_2$.

For $v_1 \in W_1^{C_1}$ and $v_2 \in {}^{A_2} \! W_2$, recall 
the definition \eqref{C2stab} of the set
$\Cbvv \subset C_2$. Similarly, let
\begin{equation}
\label{A1stab}
A_1(v_1, v_2) =\{ \al \in A_1 \mid
(v_1 c^{-1} v_2^{-1} a)^n \al \; \mbox{is defined and is in} \; A_1 \; 
\mbox{for} \; n=1, 2, \ldots\},
\end{equation}
so $A_1(v_1, v_2)$ is the largest subset of $A_1$ that is stable under
$v_1 c^{-1} v_2^{-1} a$.
Note that
\begin{equation}
\label{A1C2}
v_{2}^{-1} a \hspace{.1in} {\rm and} \hspace{.1in} 
c v_{1}^{-1} \colon \; \; 
\Aavv \lra \Cbvv
\end{equation}
are isomorphisms.
Define 
\begin{align}
\label{Kvv}
&\Kvv = \left(M_{A_1(v_1, v_2) } \times 
M_{C_2(v_1, v_2)}\right) \cap \Ad_{(e,\dot{v}_2)}^{-1} K,\\
\label{Lvv}
&\Lvv = \left(M_{A_1(v_1, v_2)} 
\times M_{C_2(v_1, v_2)}\right)
\cap \Ad_{(\dot{v}_1, e)} L,\\
\label{Qvv}
&\Qvv = \{(m, m^\prime) \in M_{\Cbvv} \times
M_{\Cbvv} \mid \exists n \in M_{A_1(v_1, v_2)} \;
\\
\nn
&\hspace{.2in}\mbox{such that} \;
(n, m) \in \Kvv \; \mbox{and} \; 
(n, m^\prime) \in \Lvv\},
\end{align}
and let $\Qvv$ act on $M_{\Cbvv}$ from the left by
\begin{equation}
\label{J1-M1}
(m, m^\prime) \cdot m_2 := m m_2 (m^\prime)^{-1}, \hspace{.2in}
(m, m^\prime) \in \Qvv, \; m_2 \in M_{\Cbvv}.
\end{equation}
 
\bth{LYmain1} \cite{LY} Let  $\A=(A_1, A_2, a, K)$ and 
$\C=(C_1, C_2, c, L)$ be 
two admissible quadruples for $G_1 \times G_2$. 
Then every $(R_\A, R_\C)$-double coset 
of $G_1 \times G_2$ is of the form
\[
R_\A(\dot{v}_1, \dot{v}_2m_2) R_\C \; \; \; \mbox{for some} \;\;
v_1 \in W_1^{C_1}, v_2 \in {}^{A_2} \! W_2, \;
m_2 \in M_{\Cbvv}.
\]
Two such double cosets $R_\A(\dot{v}_1, \dot{v}_2m_2) R_\C$ and 
$R_\A( \dot{v}'_1, \dot{v}'_2m_2^\prime) R_\C$ coincide 
if and only if $v'_i=v_i$ for $i = 1, 2$ and 
$ m_2$ and $ m^\prime_2$ 
are in the same $Q_{(v_1, v_2)}$-orbit in 
$M_{\Cbvv},$ cf. \eqref{J1-M1}.
\eth

In \cite{LY} we dealt with a slightly 
more general class of the groups
$R_\A$ for which the projections $K \to M_{A_i}$, 
for $i = 1, 2,$ do not have
to be surjective.

Recall from \eqref{partAC} that for any two admissible quadruples 
$\A$ and $\C$ for $G_1 \times G_2$, 
\[
\vvprime = R_\A (v_1, v_2M_{\Cbvv}) R_\C \subset G_1 \times G_2
\]
for $v_1 \in W_{1}^{C_1}$ and $v_2 \in {}^{A_2}\!W_2$. 
\thref{LYmain1} immediately implies:

\bco{pieces-prime} 
For any two admissible quadruples 
$\A$ and $\C$ for $G_1 \times G_2$, 
\[
G_1 \times G_2 = 
\bigsqcup_{(v_1 , v_2) \in W_{1}^{C_1} \times {}^{A_2}\!W_2}
\vvprime \hspace{.3in} (\mbox{disjoint union}).
\]
\eco

\ble{RvvT} For any $v_1 \in W_1^{C_1}$ and 
$v_2 \in {}^{A_2}\!W_2$,
\begin{equation}\label{vvAC-B}
\vvAC = R_\A (v_1, v_2 (B_2 \cap M_{\Cbvv})) R_\C,
\end{equation}
and $R_\A(v_1, v_2T_2)R_\C$ is dense in $\vvAC$, where 
$T_2 = B_2 \cap B_2^-$.
\ele

\begin{proof} Let $G_{C_2(v_1, v_2)} = M_{C_2(v_1, v_2)}/Z_{C_2(v_1, v_2)}$.
Then
\[
\sig_{(v_1, v_2)}:=\Ad_{\dot{v}_2}^{-1} 
\theta_a \Ad_{\dot{v}_1} \theta_c^{-1}
\]
defines an automorphism of $G_{\Cbvv}.$ The action of 
$Q_{(v_1, v_2)}$ on $M_{\Cbvv}$ descends to an action 
on $G_{\Cbvv}$. The orbits of the latter are exactly 
the $\sig_{(v_1,v_2)}$-twisted conjugacy classes
of $G_{\Cbvv}.$ (Here and below for a group 
$F$ and $\sig \in \Aut(F)$, by a $\sig$-twisted 
conjugacy class of $F$ we mean an
orbit of the action $g. x= g x (\sig(g))^{-1}$ of $F$ on itself.) 
\leref{RvvT} now follows from the following results in  
\cite{Mo} and \cite[Lemma 7.3]{St2}: If $G$ is a connected 
reductive algebraic group, $\sig \in \Aut(G)$,
and $B$ is a $\sig$-stable Borel subgroup 
$G$, then

(1) all $\sig$-twisted conjugacy classes of $G$ 
meet $B$;

(2) for every maximal torus $T$ of $G$ inside $B$
the union of all $\sig$-twisted conjugacy classes 
that meet $T$ is a Zariski open subset of $G$.
\end{proof}

\subsection{Stabilizer subgroups}
\label{sec-stab}
For $(g_1, g_2) 
\in G_1 \times G_2$, set
\begin{equation}
\Stab_{\A,\C}(g_1, g_2) = R_\A \cap \Ad_{(g_1, g_2)} R_\C.
\end{equation}
An explicit description of the subgroups
$\Stab_{\A,\C}(g_1, g_2)$ was given in \cite{LY}. We only recall
some facts about these groups that will be used in this paper 
(see the proof of \prref{dim-Fmi}). 

For $v_1 \in W_1^{C_1}$ and
$v_2 \in {}^{A_2} \! W_2$, let $\Abvv = a\Aavv \subset A_2.$
The following proposition was proved in \cite{LY}.

\bpr{LYmain2}  
For $q = (\dv_1, \dv_2m)$, where 
$v_1 \in W_{1}^{C_1}, v_2 \in {}^{A_2}\!W_2$,  and $m \in M_{\Cbvv},$
$\Stab_{\A, \C}(q)  \subset P_{A_1(v_1, v_2)} \times P_{A_2(v_1, v_2)}$, 
and
$\Stab_{\A, \C}(q) = \Stab_{\A, \C}^{{\rm red}}(q) 
\;\Stab_{\A, \C}^{{\rm uni}}(q)$,
where
\begin{align*}
\Stab_{\A, \C}^{{\rm red}}(q) = &\Stab_{\A, \C}(q) 
\cap (M_{\Aavv} \times M_{\Abvv}) \\
=&(M_{\Aavv} \times M_{\Abvv}) \cap K \cap 
\Ad_{(\dv_1, \dv_2m)} L  
\end{align*}
and $\Stab_{\A, \C}^{{\rm uni}}(q) = \Stab_{\A, \C}(q) 
\cap (U_{\Aavv} \times U_{\Abvv})$
has dimension equal to 
\[
\dim(U_1 \cap \Ad_{\dv_1}U_{C_1}) + \dim (U_{A_2} \cap \Ad_{\dv_2}U_2).
\]
\epr
\sectionnew{An inductive description and the geometry 
of the sets $\vvAC$}
\label{inductive-description}

In this section, we modify the inductive arguments in \cite{LY} to give 
an inductive description of the sets
$\vvAC$. Consequently, we will be able to describe the
sets $\vvAC$ as iterated bundles.
Our arguments generalize those of Lusztig in \cite{L1, L2} for 
the special case when $G_1 = G_2 = G$ and $R_\A $ is the diagonal subgroup 
of $G \times G$.
\subsection{Construction of new admissible quadruples}\label{new-construct}

For $i = 1, 2$ and $E_i, F_i \subset \Ga_i$, 
let 
\[
{}^{E_i}\!W_i^{F_i} = {}^{E_i}\!W_i \cap W_i^{F_i}.
\]
Given an admissible quadruple $\C = (C_1, C_2, c, L)$ 
for $G_1 \times G_2$, 
a subset 
$E_1$ of $\Gamma_1$, and any $y_1 \in {}^{E_1}\!W_1^{C_1}$, we 
construct another admissible quadruple $\C^{(E_1, y_1)}$ 
for $G_1 \times G_2$ as follows:
set
\begin{align*}
C_1^{(E_1, y_1)}&= E_1 \cap y_1(C_1), \; \; 
C_2^{(E_1, y_1)} = c(C_1 \cap y_1^{-1}(E_1))\\
L^{(E_1,y_1)} &= \left(M_{C_1^{(E_1, y_1)}} \times 
M_{C_2^{(E_1,y_1)}}\right) \cap 
\Ad_{(\dot{y}_1, e)} L.
\end{align*}

\ble{RE} In the above setting the following hold:

(i) The quadruple $\C^{(E_1,y_1)} :=\left(C_1^{(E_1, y_1)}, 
C_2^{(E_1, y_1)}, c y_1^{-1}, L^{(E_1, y_1)}\right)$
is an admissible quadruple for $G_1 \times G_2$.

(ii) Let
$R_{\C^{(E_1, y_1)}}$ be the subgroup of 
$G_1 \times G_2$ defined by $\C^{(E_1, y_1)}$ as in
\deref{admis}. Then 
$(P_{E_1} \times G_2) \cap \Ad_{(\dot{y}_1, e)}R_\C \subset 
R_{\C^{(E_1, y_1)}}$.
\ele

\begin{proof}
Part (i) follows directly from the definition of admissible quadruples. To
prove (ii), let
$(p, g_2) \in P_{E_1}\times G_2$ be such that
$(\dot{y}_1^{-1} p \dot{y}_1, g_2) \in R_\C$. For 
notational simplicity, set
\[
D_1 = C_1 \cap y_1^{-1}(E_1) = y_1^{-1}(C_1^{(E_1, y_1)}), \; \; \; 
D_2 = c(D_1) = C_2^{(E_1,y_1)}.
\]
Then by  \cite[(4.11) in Lemma 4.2]{LY}, 
\begin{equation}\label{pg1}
\dot{y}_1^{-1} p \dot{y}_1 \in P_{C_1} 
\cap \Ad_{\dot{y}_1}^{-1} P_{E_1}
= M_{D_1} (U_{D_1} \cap 
\Ad_{\dot{y}_1}^{-1}U_{C_1^{(E_1, y_1)} })\subset P_{D_1}.
\end{equation}
It follows from \cite[(3) in Lemma 3.4]{LY} that 
$g_2 \in P_{c(D_1)} = P_{D_2}
=P_{C_2^{(E_1, y_1)}}$, and 
\begin{align*}
(\dot{y}_1^{-1} p \dot{y}_1, g_2) \in 
(P_{D_1} \times P_{D_2}) \cap R_\C
&= ((P_{D_1} \times P_{D_2}) \cap L)(U_{C_1} \times U_{C_2})\\
&\subset ((M_{D_1} \times M_{D_2}) \cap L)(U_{D_1} \times U_{D_2}).
\end{align*}
Thus by \eqref{pg1}, 
$(\dot{y}_1^{-1} p \dot{y}_1, g_2)
\in ((M_{D_1} \times M_{D_2}) \cap L)(
\Ad_{\dot{y}_1}^{-1}U_{ C_1^{(E_1, y_1)}}  \times U_{D_2})$ and hence
$(p, g_2)
\in R_{ \C^{(E_1,y_1)}}$. 
\end{proof}
\subsection{The first step of the induction}\label{1-step}
Fix two
admissible quadruples $\A$ and $\C$ for $G_1 \times G_2$.
Given $v_1 \in W_1^{C_1}$ and $v_2 \in {}^{A_2}\!W_2$, we will use
\leref{RE} to construct a sequence of admissible quadruples for 
$G_1 \times G_2$ which will be used to give an 
inductive 
description of $\vvAC/R_\C \subset (G_1 \times G_2)/R_\C$. In this subsection, we
present the first step of the inductive description.

Consider the projection
\[
\rho_0: \; \; (G_1 \times G_2)/R_\C \lra 
(G_1 \times G_2)/(P_{C_1} \times P_{C_2}).
\]
By \cite[Proposition 8.1 and Proposition 4.1]{LY}, every
$R_\A$-orbit in 
$(G_1 \times G_2)/(P_{C_1} \times P_{C_2})$ contains exactly one
point of the form $(x_1, x_2)_\cdot (P_{C_1} \times P_{C_2})$, where
\begin{equation}\label{xx}
x_2 \in {}^{A_2}W_{2}^{C_2} \hspace{.2in} {\rm and} 
\hspace{.2in} x_1 \in {}^{a^{-1}(A_2 \cap x_2(C_2))}\!W_{1}^{C_1}.
\end{equation}
Let $\Omega_0$ be the set that consists of all pairs 
$(x_1, x_2)$ satisfying \eqref{xx}.
For $(x_1, x_2) \in \Omega_0$, let 
\begin{align*}
\Ocal^{(x_1, x_2)} &=  R_\A(x_1, x_2)_\cdot (P_{C_1} 
\times P_{C_2}) \subset (G_1 \times G_2)/(P_{C_1} \times P_{C_2}),\\
\Xcal^{(x_1, x_2)} &= \rho_{0}^{-1} (\Ocal^{(x_1, x_2)})
\subset (G_1 \times G_2)/R_\C.
\end{align*}
Then clearly,
\[
\Xcal^{(x_1, x_2)} = \{(r_1 \dx_1, \; \; r_2 \dx_2m)_\cdot 
R_\C \mid (r_1, r_2) \in R_\A, \;
m \in M_{C_2}\}.
\]
Letting $E_1 = a^{-1}\left(A_2 \cap x_2(C_2)\right)$ and $y_1 = x_1$ 
in \leref{RE}, we get
the admissible quadruple
\begin{equation}\label{newCalC}
\C^{(x_1, x_2)} :=\C^{(E_1, \, x_1)}=
(C_1^{(x_1, x_2)}, \; C_2^{(x_1, x_2)}, \; 
c^{(x_1, x_2)}, \; L^{(x_1, x_2)}) 
\end{equation}
and the corresponding subgroup 
$R_{\C^{(x_1, x_2)}}$ of $G_1 \times G_2$, where
in particular,
\begin{align}
\label{newC1}
C_1^{(x_1, x_2)} &= a^{-1}(A_2 \cap x_2(C_2)) \cap x_1(C_1), \\
\label{newC2}
C_2^{(x_1, x_2)} &= c(C_1 \cap x_{1}^{-1}a^{-1}(A_2 \cap x_2(C_2))),\\
\label{newc}
c^{(x_1, x_2)} &= cx_{1}^{-1} \colon \; \; C_1^{(x_1, x_2)} \lra 
C_2^{(x_1, x_2)}.
\end{align}

\ble{induc-1}
For any $(x_1, x_2) \in \Omega_0$, the map
$f_0 \colon \; \Xcal^{(x_1, x_2)} \to (G_1 \times G_2)/R_{\C^{(x_1, x_2)}}$
given by
\begin{equation}
\label{f0} 
f_0 \colon \; \; 
(r_1 \dx_1, \; \; r_2 \dx_2m)_\cdot R_\C \longmapsto 
(r_1, \; \; r_2 \dx_2m)_\cdot 
R_{\C^{(x_1, x_2)}}, \; \; (r_1, r_2) \in R_\A, \; m \in M_{C_2}
\end{equation}
is a well-defined $R_\A$-equivariant morphism of algebraic varieties.
\ele

\begin{proof}
Assume that $(r_1 \dx_1, r_2 \dx_2m)_\cdot R_\C =
(r_1^\prime \dx_1, r_2^\prime \dx_2m^\prime)_\cdot R_\C$
for  $(r_1, r_2), (r_1^\prime, r_2^\prime) \in R_\A$ and 
$m, m^\prime \in M_{C_2}$. Then  
\begin{equation}\label{ta} 
((r_1^\prime)^{-1}r_1, \; (r_2^\prime)^{-1}r_2) 
\in R_\A \hspace{.2in} {\rm and}
\hspace{.2in}
(\dx_{1}^{-1}(r_1^\prime)^{-1}r_1\dx_1, \; 
(m^\prime)^{-1} \dx_{2}^{-1}(r_2^\prime)^{-1}r_2\dx_2m) \in R_\C.
\end{equation}
It follows from \eqref{ta} that 
$(r_2^\prime)^{-1}r_2 
\in P_{A_2} \cap x_2(P_{C_2}) \subset P_{A_2 \cap x_2(C_2)},$
which, 
together with \eqref{ta} and \cite[(3) of Lemma 3.4]{LY}, imply 
that $(r_1^\prime)^{-1}r_1 \in P_{a^{-1}(A_2 \cap x_2(C_2))}$.
By (ii) of \leref{RE}, 
$((r_1^\prime)^{-1}r_1, \; 
(m^\prime)^{-1} \dx_{2}^{-1}(r_2^\prime)^{-1}r_2\dx_2m) 
\in R_{\C^{(x_1, x_2)}}$.
Thus $f_0$ is well-defined. 

Clearly $f_0$ is $R_\A$-equivariant. To see that
$f_0$ is a morphism,  
note that since the fiber of $\rho_0$
over the point $(\dot{x}_1, \dot{x}_2)_\cdot (P_{C_1} \times P_{C_2})$
is $(\dot{x}_1, \dot{x}_2)(P_{C_1} \times P_{C_2})/R_\C$, 
\[
\Xcal^{(x_1, x_2)} \cong R_\A 
\times_{R_1}(\dot{x}_1, \dot{x}_2)(P_{C_1} \times P_{C_2})/R_\C,
\]
where 
$R_1 = R_\A \cap \Ad_{(\dot{x}_1, \dot{x}_2)}(P_{C_1} \times P_{C_2})$
is the stabilizer subgroup of $R_\A$ at the point
$(\dot{x}_1, \dot{x}_2)_\cdot (P_{C_1} \times P_{C_2})$. Let
$Y_2 = \{m \in M_{C_2} \mid (e, m) \in L\} \subset Z_{C_2}$.
Then the inclusion map of $\{e\} \times M_{C_2} \hookrightarrow P_{C_1}
\times P_{C_2}$ induces an isomorphism
\[
\phi_1: \; \; (\dx_1, \dx_2)(\{e\} \times M_{C_2})/(\{e\} \times Y_2) 
\stackrel{\cong}{\longrightarrow}
(\dx_1, \dx_2)(P_{C_1} \times P_{C_2})/R_\C.
\]
Note that $\{e\} \times Y_2 \subset L^{(x_1, x_2)}\subset R_{\C^{(x_1, x_2)}}$
so we have the projection 
\[
\phi_2: \; \; (G_1 \times G_2)/(\{e\} \times Y_2) \to 
(G_1 \times G_2)/R_{\C^{(x_1, x_2)}}.
\]
Consider the morphism
\begin{align*}
\phi_3: & R_\A \times (G_1 \times G_2)/(\{e\} \times Y_2) \longrightarrow 
(G_1 \times G_2)/(\{e\} \times Y_2),\\
& ((r_1, r_2), (g_1, g_2)(\{e\} \times Y_2)) \longmapsto (r_1 \dot{x}_1^{-1} g_1, \, 
r_2g_2)(\{e\} \times Y_2).
\end{align*}
The composition of $\phi_3$ with $\phi_2$ gives rises to a morphism
\[
R_\A \times (\dot{x}_1, \dot{x}_2)(\{e\} \times M_{C_2})/(\{e\} \times Y_2)
\longrightarrow (G_1 \times G_2)/R_{\C^{(x_1, x_2)}}.
\]
Let $R_1$ act on $(\dot{x}_1, \dot{x}_2)(\{e\} \times M_{C_2})/(\{e\} \times Y_2)$
so that the isomorphism $\phi_1$ is $R_1$-equivariant. Then the
well-definedness of $f_0$ implies that 
\[
f_0: \; \; \chi^{(x_1, x_2)} \cong 
R_\A \times_{R_1} (\dot{x}_1, \dot{x}_2)(\{e\} \times M_{C_2})/(\{e\} \times Y_2)
\longrightarrow (G_1 \times G_2)/R_{\C^{(x_1, x_2)}}
\]
is a morphism of varieties.
\end{proof}

Let now $v_1 \in W_1^{C_1}$ and $v_2 \in {}^{A_2}\!W_2$. By 
\cite[Proposition 2.7.5]{Car} (see \leref{car-prop275}
in the Appendix), 
$v_2$ can be uniquely written as
\begin{equation}\label{v2}
v_2 = x_2 u_2, \hspace{.2in} \mbox{where} \hspace{.2in}
x_2 \in {}^{A_2}\!W_2^{C_2}, \; \; \mbox{and} \; \; 
u_2 \in {}^{C_2 \cap x_2^{-1}(A_2)}\!W_{C_2}.
\end{equation}
Given the decomposition $v_2 = x_2u_2$, write $v_1$ as 
\begin{equation}
\label{v1}
v_1 =u_1 x_1,  \hspace{.2in} \mbox{where} \hspace{.2in} 
x_1 \in {}^{a^{-1}(A_2 \cap x_2(C_2))}\!W_{1}^{C_1}, \; \; 
\mbox{and} \; \; u_1 \in 
W_{a^{-1}(A_2 \cap x_2(C_2))}^{C_{1}^{(x_1, x_2)}}
\end{equation}
with $C_{1}^{(x_1, x_2)}$ given in \eqref{newC1}.
Here and below, for $i = 1, 2$ and $E_i\subset D_i \subset
\Ga_i$, we let 
\[
W_{D_i}^{E_i} = W_{D_i} \cap W_i^{E_i}.
\]
It is easy to see that $\vvAC/R_\C \subset \chi^{(x_1, x_2)}$.
Since $u_1 \in 
W_{a^{-1}(A_2 \cap x_2(C_2))}^{ C_1^{(x_1, x_2)}} 
\subset W_{1}^{ C_1^{(x_1, x_2)}}$, 
we have the set
\[
[u_1, v_2]_{\A,  \C^{(x_1, x_2)}} / R_{ \C^{(x_1, x_2)}} =
R_\A \left(u_1, v_2 M_{ C_2^{(x_1, x_2)}(u_1, v_2)} \right)
\cdot 
R_{ \C^{(x_1, x_2)}} \subset (G_1 \times G_2)/R_{ \C^{(x_1, x_2)}},
\]
where $ C_2^{(x_1, x_2)}(u_1, v_2)$ is the largest subset of 
$C_2^{(x_1, x_2)}$ that
is stable under the map 
\[
v_{2}^{-1}a u_1 ( c^{(x_1, x_2)})^{-1}  = 
v_{2}^{-1}a v_1 c^{-1}.
\]
An argument similar to that in the proof of \cite[Lemma 5.3]{LY} 
shows that $ C_2^{(x_1, x_2)}(u_1, v_2) = \Cbvv$.
The following \prref{zx} now follows directly from \thref{LYmain1}.

\bpr{zx}
Let   
$v_1 \in W_1^{C_1}$ and $v_2 \in {}^{A_2}\!W_2$ 
have the decompositions \eqref{v1} and \eqref{v2}. 
Then $\vvAC$ consists of those $(g_1, g_2) \in G_1 \times G_2$ 
for which 
\[
(g_1, g_2) . R_\C \in \Xcal^{(x_1, x_2)} \hspace{.2in} {\rm and}
\hspace{.2in} f_0((g_1, g_2) . R_\C) \in 
[u_1, v_2]_{\A,  \C^{(x_1, x_2)}}/R_{\C^{(x_1, x_2)}}.
\]
Moreover, 
$f_0 \colon  [v_1, v_2]_{\A, \C}/R_\C 
\to [u_1, v_2]_{\A,  \C^{(x_1, x_2)}}/R_{\C^{(x_1, x_2)}}$
induces a bijection between the $R_\A$-orbits in $\vvprime/R_\C$ 
and the $R_\A$-orbits in 
$[u_1, v_2]_{\A,  \C^{(x_1, x_2)}}/R_{\C^{(x_1, x_2)}}$.
\epr
\subsection{Inductive description of the sets $\vvAC/R_\C$}\label{induc-steps}

Given $v_1 \in W_1^{C_1}$ and $v_2 \in {}^{A_2}\!W_2$, by repeating the 
construction in \S \ref{1-step}, we get a sequence
\[
\C^{(i)} = (C_{1}^{(i)},  C_{2}^{(i)},  c^{(i)}, 
L^{(i)}), \; \; \; i \geq 0,
\]
of admissible quadruples for $G_1 \times G_2$ and a sequence
$u_1^{(i)} \in W_1^{C_1^{(i)}}$, $i \geq 0, \dots$
which gives rise to the sequence of double cosets
\[
\Zcal^{(i)} = [u_{1}^{(i)}, \; v_2]_{\A, \C^{(i)}}/R_{\C^{(i)}} 
\subset (G_1 \times G_2)/R_{\C^{(i)}},
\hspace{.2in}
i \geq 0.
\]
Here
$\C^{(0)} = \C$, $u_1^{(0)} = v_1$, $\Zcal^{(0)} = \vvAC/R_\C$,
$\C^{(1)} = \C^{(x_1, x_2)}$ as in 
\eqref{newCalC}, and $u_1^{(1)} = u_1$ as in \eqref{v1}. In general,
once $(\C^{(i)}, u_1^{(i)})$ is given,
$(\C^{(i+1)}, u_1^{(i+1)})$ is constructed from $(\C^{(i)}, u_1^{(i)})$
in the same way as $(\C^{(x_1, x_2)}, u_1)$ 
was from $(\C, v_1)$ in \S \ref{1-step}. 
Namely,
first decompose $v_2$ as the unique product
$v_2 = x_{2}^{(i)}u^{(i)}_{2}$ with
\begin{equation}\label{xu}
x_{2}^{(i)} \in {}^{A_2}\!W_{2}^{C_{2}^{(i)}}, \; \; u^{(i)}_{2} \in 
{}^{C_{2}^{(i)}\cap (x_{2}^{(i)})^{-1}(A_2) }\!W_{C_{2}^{(i)}}.
\end{equation}
Then decompose $u_{1}^{(i)}$ as the unique product
$u_{1}^{(i)} = u_{1}^{(i+1)} x_{1}^{(i)},$
where
\[
x_{1}^{(i)} \in {}^{a^{-1}(A_2 \cap 
x_{2}^{(i)}(C_{2}^{(i)}))}\!W_1^{C_{1}^{(i)}},
\; \; \; u_1^{(i+1)} \in W_{a^{-1}(A_2 \cap x_{2}^{(i)}
(C_{2}^{(i)}))}^{a^{-1}(A_2 \cap x_{2}^{(i)}(C_{2}^{(i)})) 
\cap x_{1}^{(i)} (C_{1}^{(i)})}.
\]
The admissible quadruple $\C^{(i+1)}$ is constructed as in 
\S \ref{new-construct}
by taking  $\C$ to be $\C^{(i)}$, $E_1$ to be 
$a^{-1}(A_2 \cap x_{2}^{(i)}(C_{2}^{(i)}))$ 
and $y_1$ to be $x_{1}^{(i)}$.

For $i \geq 0$, 
let $\rho_i$ be the natural projection
\[
\rho_i \colon \; \; (G_1 \times G_2)/R_{\C^{(i)}} \lra 
(G_1 \times G_2)/(P_{C^{(i)}_{1}} \times P_{C^{(i)}_{2}}),
\]
$\Ocal^{(i)} :=  R_\A (x_{1}^{(i)}, x_{2}^{(i)})_\cdot 
(P_{C^{(i)}_{1}} \times P_{C^{(i)}_{2}})$, and
\[
\Xcal^{(i)}  := \rho_i^{-1}(\Ocal^{(i)}) \subset
(G_1 \times G_2)/R_{\C^{(i)}}.
\]
Then $\Xcal^{(i)}$ is a locally closed subset of 
$(G_1 \times G_2)/R_{\C^{(i)}}$. 
By \leref{induc-1} and \prref{zx}, we have a well-defined 
$R_\A$-equivariant morphism
$f_i \colon \; \; \Xcal^{(i)} \lra (G_1 \times G_2)/R_{\C^{(i+1)}}$
such that
\begin{equation}\label{Zi}
\Zcal^{(i)} \subset \Xcal^{(i)} \hspace{.2in }{\rm and} \hspace{.2in} 
\Zcal^{(i)} = f_{i}^{-1}\left(\Zcal^{(i+1)}\right), 
\hspace{.2in} \forall i \geq 0.
\end{equation}
Let
$i_0 \geq 0$ be the smallest integer such that 
$C_{2}^{(i_0+1)} = C_{2}^{(i_0)}$.
It is then easy
to see that 
\begin{equation}\label{Ci}
\C^{(i_0)} =
\C^{(\infty)} := 
(A_1(v_1, v_2), \, C_2(v_1, v_2), \, cv_{1}^{-1}, \, \Lvv) \; \; 
\mbox{and}  \; \; u_{1}^{(i)} = e, 
\hspace{.1in} \forall i \geq i_0 + 1,
\end{equation}
where
$L_{(v_1, v_2)}$ is given by \eqref{Lvv}. 
Set $\Zcal^{(\infty)}=\Zcal^{(i_0+1)}$ and
$\Xcal^{(\infty)} = \Xcal^{(i_0+1)}$. 
It follows from \eqref{Ci} that 
\[
\Zcal^{(i)} = \Zcal^{(\infty)} = \Xcal^{(\infty)} = 
[e, v_2]_{\A, \C^{(\infty)}}/R_{\C^{(\infty)}}, \hspace{.2in}
\forall i \geq i_0 + 1.
\] 

\bpr{Zi-locallyclosed}
For every $i \geq 0$, $\Zcal^{(i)}$ is locally closed in
$(G_1 \times G_2)/R_{\C^{(i)}}$. 
\epr

\begin{proof}  By the inductive description of the sets 
$\Zcal^{(i)}$ in \eqref{Zi}, if  for some $i \geq 0$, 
$\Zcal^{(i+1)}$ is locally closed in 
$(G_1 \times G_2)/R_{\C^{(i+1)}}$, then 
$\Zcal^{(i)}$ is locally closed in
$(G_1 \times G_2)/R_{\C^{(i)}}$ because $f_i$ is
a morphism. Since $\Zcal^{(\infty)}=
\Xcal^{(\infty)}$ is locally closed in 
$(G_1 \times G_2)/R_{\C^{(\infty)}}$, it follows that 
$\Zcal^{(i)}$ is locally closed in
$(G_1 \times G_2)/R_{\C^{(i)}}$ for every $i \geq 0$.
\end{proof}

Let $\rho_{\infty} := \rho_{i_0+1} \colon \; 
(G_1 \times G_2)/R_{\C^{(\infty)}} \to
(G_1 \times G_2)/(P_{\Aavv} \times P_{\Cbvv})$, and
\[
\Ocal^{(\infty)} :=\Ocal^{(i_0+1)}=R_\A(e, v_2)_\cdot 
(P_{\Aavv} \times P_{\Cbvv}),
\]
so
that $\Zcal^{(\infty)} = \Xcal^{(\infty)}=
\rho_{\infty}^{-1}(\Ocal^{(\infty)})$. Then we have the projection
$\rho_{\infty} \colon 
\Zcal^{(\infty)} \to \Ocal^{(\infty)}.$
Note that the stabilizer subgroup of $R_\A$ at the point
\[
(e, v_2)_\cdot (P_{\Aavv} \times P_{\Cbvv}) \in 
(G_1 \times G_2)/(P_{\Aavv} \times P_{\Cbvv})
\]
is
\begin{equation}\label{Rinclu}
R_\A \cap \left(P_{\Aavv} \times \Ad_{\dv_2}P_{\Cbvv}\right) \subset
R_\A \cap (P_{\Aavv} \times P_{\Abvv}),
\end{equation}
where $A_2(v_1, v_2) = a A_1(v_1, v_2)$. For notational simplicity, set
\begin{equation}\label{Rvv}
\Rvv = R_\A \cap (P_{\Aavv} \times P_{\Abvv}).
\end{equation}
Let 
$\rho_{\infty}^{\prime} \colon  \Ocal^{(\infty)} \to
R_\A / \Rvv$ be the projection 
induced by the inclusion map in \eqref{Rinclu}. 
Then we have the sequence of $R_\A$-equivariant
morphisms
\begin{equation}
\vvprime/R_\C = \Zcal^{(0)}
\stackrel{f_0}{\lra} \Zcal^{(1)}
\stackrel{f_1}{\lra} \cdots \stackrel{f_{i_0-1}}{\lra} 
\Zcal^{(i_0)}
\stackrel{f_{i_0}}{\lra}  \Zcal^{(\infty)} \\
\stackrel{\rho_\infty}{\lra} \Ocal^{(\infty)}
\stackrel{\rho_\infty^\prime}{\lra}R_\A/R_{(v_1, v_2)}.
\label{fibration-sequence}
\end{equation}

\ble{Zinfty2} The quotient
$R_\A/\Rvv$ is isomorphic to the 
flag variety $M_{A_1}/(M_{A_1} \cap P_{\Aavv})$
of $M_{A_1}$. The fibers of $\rho_\infty:
\Zcal^{(\infty)} \to \Ocal^{(\infty)}$
are isomorphic to $M_{\Cbvv}/Y_2$, 
where $Y_2 =\{m \in M_{C_2} \mid (e, m) \in L\} \subset Z_{C_2}$,
and the fibers of 
$\rho_\infty^\prime: \Ocal^{(\infty)}
\to R_\A/R_{(v_1, v_2)}$ are isomorphic to 
$U_{2} \cap \Ad_{\dv_2}U_{2}^{-}$.
\ele

\begin{proof} Consider the group homomorphism
\[
p \colon \; \; R_\A \lra M_{A_1}\colon \; \; (r_1, r_2) \longmapsto
m_1 \hspace{.2in} {\rm if} \hspace{.1in} r_1 = m_1 u_1 \hspace{.1in}
{\rm for} \hspace{.1in} m_1 \in M_{A_1}, u_1 \in U_{A_1}.
\]
Since $p$ is surjective, the action of $R_\A$ on 
$M_{A_1}/(M_{A_1} \cap P_{\Aavv})$  through the homomorphism $p$
is transitive. The stabilizer subgroup of $R_\A$ at the 
point $e_\cdot (M_{A_1} \cap P_{\Aavv})$ is $\Rvv$ by \cite[(3.7) in Lemma 3.4]{LY}. 
Thus $R_\A/\Rvv$ is isomorphic to 
$M_{A_1}/(M_{A_1} \cap P_{\Aavv})$. 

The fibers of $\rho_\infty: \Zcal^{(\infty)} \to \Ocal^{(\infty)} $ 
are clearly isomorphic to 
\[
(P_{A_1(v_1, v_2)} \times P_{C_2(v_1, v_2)})/R_{\C^{(\infty)}} \cong 
(M_{A_1(v_1, v_2)} \times M_{C_2(v_1, v_2)})/L_{(v_1, v_2)}.
\]
Note that $\{e\} \times M_{C_2(v_1, v_2)}$ acts on  
$(M_{A_1(v_1, v_2)} \times M_{C_2(v_1, v_2)})/L_{(v_1, v_2)}$
transitively, and
\[
(\{e\} \times M_{C_2(v_1, v_2)}) \cap L_{(v_1, v_2)} \cong \{e\} \times Y_2.
\]
Thus the fibers of $\rho_\infty: \Zcal^{(\infty)} \to \Ocal^{(\infty)}$
are isomorphic to $M_{\Cbvv}/Y_2$.
The fibers of $\rho_\infty^\prime$ are isomorphic to 
\[
R_{(v_1, v_2)}/R_\A \cap (P_{A_1(v_1, v_2)} \times \Ad_{\dv_2} 
P_{C_2(v_1, v_2)}) \cong U_2 \cap \Ad_{\dv_2}U_2^-.
\] 
\end{proof}
\subsection{The set $\vvAC/R_\C$ as an iterated 
fiber bundle}\label{iterated}
 
Assume the setting from \S \ref{induc-steps}. 
We will show in this subsection
that each morphism $f_i: \Zcal^{(i)} \to \Zcal^{(i+1)}$
in \eqref{fibration-sequence} has fibers isomorphic to an affine space. 

Fix $m \in M_{C_2(v_1, v_2)}$. 
For $i \geq 0$, let $z_m^{(i)} = 
(\dot{u}_1^{(i)}, \dot{v}_2 m)_\cdot R_{\C^{(i)}}
\in \Zcal^{(i)}$ and
\[
S_m^{(i)} = R_\A \cap \Ad_{(\dot{u}_1^{(i)}, \dot{v}_2m)} R_{\C^{(i)}}
\]
be the stabilizer subgroup of $R_\A$ at $z_m^{(i)} \in \Zcal^{(i)}$.
Set
\[
F_{m}^{(i)} = f_i^{-1}(z_m^{(i+1)}) \subset \Zcal^{(i)}.
\]

\ble{Fmi}
For any $m \in M_{C_2(v_1, v_2)}$ and $i \geq 0$, 
$S_m^{(i)} \subset S_m^{(i+1)}$, 
\begin{equation}\label{eq-Fmi}
F_{m}^{(i)} = {S_m^{(i+1)}}_\cdot z_m^{(i)} 
\cong S_m^{(i+1)}/S_m^{(i)}.
\end{equation}
\ele

\begin{proof} Let $i \geq 0$. 
Since $f_i(z_m^{(i)}) = z_m^{(i+1)}$ and since $f_i$ is 
$R_\A$-equivariant, we 
have $S_m^{(i)} \subset S_m^{(i+1)}$,\ 
and  ${S_m^{(i+1)}}_\cdot z_m^{(i)}
\subset F_{m}^{(i)}$. 
It suffices to prove \eqref{eq-Fmi} for $i = 0$, 
and we only need to show that
$F_{m}^{(0)} \subset {S_m^{(1)}}_\cdot z_m^{(0)}$.

Recall that $\C^{(1)} = \C^{(x_1, x_2)}$ as in 
\eqref{newCalC}, and $u^{(1)} = u_1$
as in the $v_1 = u_1x_1$ decomposition in \eqref{v1}. 
To show that $F_{m}^{(0)} \subset {S_m^{(1)}}_\cdot z_m^{(0)}$, 
assume that
$(r_1, r_2) \in R_\A$ and $m^\prime \in M_{C_2(v_1, v_2)}$ 
are such that
$(r_1 \dot{v}_1, r_2 \dot{v}_2 m^\prime)_\cdot R_\C \in F_m^{(0)}$. 
Then
\[
R_\A(\dot{u}_1, \dot{v}_2 m^\prime)_\cdot R_{\C^{(x_1, x_2)}} =
R_\A(\dot{u}_1, \dot{v}_2 m)_\cdot R_{\C^{(x_1, x_2)}}.
\]
Applying \thref{LYmain1} to the quadruples 
$(\A, \C^{(x_1, x_2)})$ and
noting that the largest subset of $A_1$ that is stable under
$u_1 (c^{(x_1, x_2)})^{-1}v_2^{-1}a = v_1 c^{-1} v_{2}^{-1} a$ is 
$\Aavv$ and 
that
\[
(M_{\Aavv} \times M_{\Cbvv}) \cap \Ad_{(\du_1, e)} L^{(x_1, x_2)} =
L_{(v_1, v_2)},
\]
we know that there exist $m_1 \in M_{\Aavv}$ and $n, n^\prime \in M_{\Cbvv}$
with $(m_1, n) \in K_{(v_1, v_2)}$ and 
$(m_1, n^\prime) \in L_{(v_1, v_2)}$ such that
$m^\prime = n m (n^\prime)^{-1}$, where $\Kvv$ and $\Lvv$ are 
respectively given by \eqref{Kvv} and \eqref{Lvv}. Thus
\begin{align*}
(r_1 \dot{v}_1, r_2 \dot{v}_2 m^\prime)_\cdot R_\C
&= (r_1 \dv_1,\; r_2 \dv_2 n m (n^\prime)^{-1})_\cdot R_\C
= (r_1 m_1 \dv_1, \; r_2 \dv_2 nm)_\cdot R_\C\\
&= (r_1, r_2)(m_1, \Ad_{\dv_2}(n)) (\dv_1, \;\dv_2 m)_\cdot R_\C=
(r_1^\prime, r_2^\prime) z_m^{(0)},
\end{align*}
where 
$(r_1^\prime, r_2^\prime) = (r_1, r_2)(m_1, \Ad_{\dv_2}(n)) \in R_\A$.
On the other hand, it follows from  
$(r_1 \dot{v}_1, r_2 \dot{v}_2 m^\prime)_\cdot R_\C \in F_m^{(0)}$ that
\[
z_m^{(1)} = 
(r_1 \dot{u}_1,\; r_2\dv_2 m^\prime)_\cdot R_{\C^{(x_1, x_2)}}
=(r_1\dot{u}_1,\; r_2 \dv_2 n m (n^\prime)^{-1})_\cdot 
R_{\C^{(x_1, x_2)}}=
(r_1^\prime, r_2^\prime)z_{m}^{(1)},
\]
so $(r_1^\prime, r_2^\prime)  \in S_{m}^{(1)}$, and thus
$(r_1 \dot{v}_1, r_2 \dot{v}_2 m^\prime)_\cdot R_\C 
\in {S_{m}^{(1)}}_\cdot z_m^{(0)}$. Hence
$F_{m}^{(0)} \subset {S_m^{(1)}}_\cdot z_m^{(0)}$.
\end{proof}

\bpr{dim-Fmi}
For each $i \geq 0$, the fibers of the morphism 
$f_i: \Zcal^{(i)} \to \Zcal^{(i+1)}$ 
are
isomorphic to the affine space of
dimension
\begin{equation}\label{dim}
\dim\left(U_1 \cap \Ad_{\dot{u}_{1}^{(i+1)}} U_{C_1^{(i+1)}}\right)-
\dim\left(U_1 \cap \Ad_{\dot{u}_{1}^{(i)}} U_{C_1^{(i)}}\right).
\end{equation}
\epr

\begin{proof}
Let $i \geq 0$ and $m \in M_{C_2(v_2, v_2)}$. 
By \prref{LYmain2}, we have the semi-direct product decompositions
$S_m^{(i)} = S_m^{(i), {\rm red}}S_m^{(i), {\rm uni}}$, where
\[
S_m^{(i), {\rm red}} = S_m^{(i)} \cap (M_{A_1(v_1, v_2)} 
\times M_{A_2(v_1, v_2)}), 
\hspace{.2in}
S_m^{(i), {\rm uni}} = S_m^{(i)} \cap (U_{A_1(v_1, v_2)} 
\times U_{A_2(v_1, v_2)}).
\]
It also follows from \prref{LYmain2} that 
$S_m^{(i), {\rm red}}=S_m^{(i+1), {\rm red}}$.
Thus by \eqref{eq-Fmi} and by \prref{LYmain2}, 
$F_m^{(i)}\cong S_m^{(i+1), {\rm uni}} /S_m^{(i), {\rm uni}}$ 
is isomorphic to an affine space of 
dimension \eqref{dim}.
\end{proof}

For $v_1 \in W_1^{C_1}$ and $v_2 \in {}^{A_2}\!W_2$, 
consider the composition 
\[
p_0 = \rho_\infty^\prime \circ \rho_\infty \circ f_{i_0} \circ f_{i_0-1} 
\circ \cdots \circ f_1 \circ f_0: \; \vvAC \lra R_\A/R_{(v_1, v_2)}.
\]
Since for $(r_1, r_2) \in R_\A$ and $m \in M_{C_2(v_1, v_2)}$,
\[
p_0((r_1, r_2)(\dv_1, \dv_2m)_\cdot R_\C) = (r_1, r_2)_\cdot R_{(v_1, v_2)},
\]
it follows that the fiber of 
$p_0$ over the point $e_\cdot R_{(v_1, v_2)} \in R_\A/R_{(v_1, v_2)}$ is
\[
X_{(v_1, v_2)} \stackrel{{\rm def}}{=}
R_{(v_1, v_2)}(v_1, v_2 M_{\Cbvv})_\cdot R_\C \subset \vvAC/R_\C.
\]
Introduce
\[
P^\prime_{A_2(v_1, v_2)} = M_{A_2(v_1, v_2)} (U_{A_2(v_1, v_2)} 
\cap \Ad_{\dv_2}U_2^-) \subset P_{A_2(v_1, v_2)}.
\]
Since $M_{A_2(v_1, v_2)}$ normalizes 
$U_{A_2(v_1, v_2)} \cap \Ad_{\dv_2}U_2^-$, $P^\prime_{A_2(v_1, v_2)}$ is a subgroup of
$P_{A_2(v_1, v_2)}$.

\bco{properties}
For any $v_1 \in W_1^{C_1}$ and $v_2 \in {}^{A_2}\!W_2$, (i) $X_{(v_1, v_2)}$
is isomorphic to  
\begin{align*}
& (P_{A_1(v_1, v_2)} \times P_{A_2(v_1, v_2)}) 
(\dv_1, \dv_2)_\cdot R_\C
=(U_{A_1(v_1, v_2)} \times P^\prime_{A_2(v_1, v_2)}) (\dv_1, \dv_2)_\cdot R_\C\\
&\cong \left(U_{\Aavv}/(U_{\Aavv} \cap \Ad_{\dv_1}U_{C_1})\right) 
\times (U_2 \cap \Ad_{\dv_2}U_2^-) \times 
(M_{\Cbvv}/Y_2).
\end{align*}
In particular, $X_{(v_1, v_2)}$ is a smooth locally closed subset of $(G_1 \times G_2)/R_\C$;

(ii) the action map of $R_\A$ on $\vvAC$ gives rise to an
isomorphism
\[
R_\A \times_{\Rvv} X_{(v_1, v_2)} \lra \vvAC/R_\C.
\]
In particular, 
$\vvAC/R_\C$ is a smooth irreducible locally closed subset of 
$(G_1 \times G_2)/R_\C$.
\eco

Since $G_1 \times G_2 \to (G_1 \times G_2)/R_\C$ is 
a locally trivial fibration, $\vvAC$ is a smooth irreducible locally 
closed subset of $G_1 \times G_2$. Note also that 
$U_{\Aavv} \cap \Ad_{\dv_1}U_{C_1} = U_1 \cap \Ad_{\dv_1}U_{C_1}$, so (i) of
\thref{main1} follows now from \coref{properties}.

\begin{proof} (i) Clearly $X_{(v_1, v_2)} \subset 
(P_{A_1(v_1, v_2)} \times P_{A_2(v_1, v_2)}) (\dv_1, \dv_2)_\cdot R_\C$.
Since 
\begin{align*}
P_{A_1(v_1, v_2)} \times P_{A_2(v_1, v_2)}&= R_{(v_1, v_2)} 
(\{e\} \times P_{\Abvv} \cap M_{A_2})\\&= R_{(v_1, v_2)}
(\{e\} \times M_{\Abvv}(U_{\Abvv}  \cap M_{A_2})),
\end{align*} 
we have
\[
(P_{A_1(v_1, v_2)} \times P_{A_2(v_1, v_2)}) (\dv_1, \dv_2)_\cdot R_\C
=R_{(v_1 v_2)} (\dv_1, \dv_2M_{C_2(v_1, v_2)} 
\Ad_{\dv_2}^{-1} (U_{\Abvv} \cap M_{A_2}))_\cdot R_\C.
\]
Since $v_2 \in {}^{A_2}\!W_2$, $\Ad_{\dv_2}^{-1} (U_{\Abvv} \cap M_{A_2})
\subset U_2$. Thus
\[
(P_{A_1(v_1, v_2)} \times P_{A_2(v_1, v_2)}) (\dv_1, \dv_2)_\cdot R_\C
\subset R_{(v_1 v_2)} (\dv_1, \dv_2P_{C_2(v_1, v_2)})_\cdot R_\C.
\]
By following the definition of $f_i$ for each $i \geq 0$, one sees that 
$(\dv_1, \dv_2P_{C_2(v_1, v_2)})_\cdot R_\C \subset X_{(v_1, v_2)}.$
Thus $(P_{A_1(v_1, v_2)} \times P_{A_2(v_1, v_2)}) (\dv_1, \dv_2)_\cdot R_\C
=X_{(v_1, v_2)}$.
Now since
\begin{align*}
P_{A_2(v_1, v_2)} &=M_{A_2(v_1, v_2)} (U_{A_2(v_1, v_2)} \cap \Ad_{\dv_2}U_2^-)
(U_{A_2(v_1, v_2)} \cap \Ad_{\dv_2}U_2) \\
&= P^\prime_{A_2(v_1, v_2)} 
(U_{A_2(v_1, v_2)} \cap \Ad_{\dv_2}U_2),
\end{align*}
we have $P_{A_2(v_1, v_2)}\dv_2= P^\prime_{A_2(v_1, v_2)}\dv_2 
(U_2 \cap \Ad_{\dv_2}^{-1}U_{A_2(v_1, v_2)})$.
Thus 
\[
X_{(v_1, v_2)} = (P_{A_1(v_1, v_2)}\dv_1, \; P^\prime_{A_2(v_1, v_2)} \dv_2
(U_2 \cap \Ad_{\dv_2}^{-1}U_{A_2(v_1, v_2)} \cap M_{C_2})_\cdot R_\C.
\]
Using the facts that $v_2^{-1} A_2(v_1, v_2) = C_2(v_1, v_2)$ and 
$v_1c^{-1}C_2(v_1, v_2) = A_1(v_1, v_2)$, we have
\[
X_{(v_1, v_2)} = (P_{A_1(v_1, v_2)}\dv_1, \; P^\prime_{A_2(v_1, v_2)} 
\dv_2)_\cdot R_\C.
\]
Using again the fact that $cv_1^{-1}A_1(v_1, v_2) = C_2(v_1, v_2)$, we get
\[
X_{(v_1, v_2)} = (U_{A_1(v_1, v_2)} \times P^\prime_{A_2(v_1, v_2)}) 
(\dv_1, \dv_2)_\cdot R_\C.
\]
Let $S = (U_{\Aavv} \times P^\prime_{\Abvv}) \cap \Ad_{(\dv_1, \dv_2)} R_\C$.
We now determine $S$. Suppose that $u \in U_{\Aavv}$ and $p^\prime \in P^\prime_{\Abvv}$
are such that $(u, p^\prime) \in S$. Then
$(\dv_1^{-1} u \dv_1, \dv^{-1}_2 p^\prime \dv_2) \in R_\C$. Hence
\[
\dv_1^{-1} u \dv_1 \in (\Ad_{\dv_1}^{-1}U_{\Aavv}) \cap P_{C_1} = 
(\Ad_{\dv_1}^{-1}U_{\Aavv}) \cap U_{C_1}.
\]
Thus $\dv^{-1}_2 p^\prime \dv_2 \in Y_2 U_{C_2}$. Since $\dv^{-1}_2 p^\prime \dv_2
\in M_{\Cbvv} U_2^- \subset P_{\Cbvv}^-$, we must have
$\dv^{-1}_2 p^\prime \dv_2 \in Y_2$. Thus $(u, p^\prime) \in 
(U_{\Aavv} \cap \Ad_{\dv_1} U_{C_1} ) \times \Ad_{\dv_2}Y_2$. 
Since
\[
(U_{\Aavv} \cap \Ad_{\dv_1} U_{C_1} ) \times \Ad_{\dv_2}Y_2 \subset S,
\]
we have $S = (U_{\Aavv} \cap \Ad_{\dv_1} U_{C_1} ) \times \Ad_{\dv_2}Y_2$. Thus
$X_{(v_1, v_2)} $ is isomorphic to
\begin{align*}
&\!(U_{A_1(v_1, v_2)} \times P^\prime_{A_2(v_1, v_2)})/S \\
&\cong \!\!\left(U_{\Aavv}/(U_{\Aavv} \!\cap \!\Ad_{\dv_1}U_{C_1})\right) \!\times \!
(U_{\Abvv} \cap \Ad_{\dv_2}U_2^-) \!\times \!(M_{\Abvv}/\Ad_{\dv_2}Y_2) \\
& \cong \!\!\left(U_{\Aavv}/(U_{\Aavv} \cap \Ad_{\dv_1}U_{C_1})\right) \times 
(U_2 \cap \Ad_{\dv_2}U_2^-) \times (M_{\Cbvv}/Y_2).
\end{align*}
 
(ii) Since $p_0$ is $R_\A$-equivariant, by \cite[Lemma 4]{Sl}, 
the action map of $R_\A$ on $\vvAC$
induces an isomorphism
\[
R_\A \times_{\Rvv} X(v_1, v_2) \lra \vvAC
\]
of $R_\A$-varieties. 
\end{proof}

\subsection{Another description of the strata $\vvAC/R_\C$}\label{vvAC-BB}
\bpr{2dc-plus-2} 
For any $v_1 \in W_1^{C_1}, \; v_2 \in {}^{A_2} \! W_2$, one has
\begin{equation}\label{vvBB}
[v_1, v_2]_{\A, \C}/R_\C = R_\A (B_1\times B_2)(v_1, v_2)_\cdot R_\C.
\end{equation}
\epr

\begin{proof} By \leref{RvvT}, we have
\begin{multline*}
\vvAC/R_\C = R_\A (v_1, v_2 (B_2 \cap M_{\Cbvv}))_\cdot R_\C 
\subset \\
R_\A (v_1, v_2 B_2)_\cdot R_\C 
= R_\A (B_1\times B_2)(v_1, v_2)_\cdot R_\C.
\end{multline*}
It remains to show that 
$R_\A (B_1\times B_2)(v_1, v_2)_\cdot R_\C \subset \vvAC/R_\C.$
Write $v_2 = x_2u_2$ and $v_1 = u_1x_1$ as in 
\eqref{v2} and \eqref{v1} and let 
$f_0 \colon \; \Xcal^{(x_1, x_2)} \to 
(G_1 \times G_2)/R_{\C^{(x_1, x_2)}}$
be as in \leref{induc-1}. Then 
\begin{align*}
R_\A (B_1\times B_2)(v_1, v_2)_\cdot R_\C &=
R_\A(v_1, v_2 (B_2 \cap M_{C_2}))_\cdot R_\C \subset 
\Xcal^{(x_1, x_2)},\\
f_0(R_\A (B_1 \times B_2) (v_1, v_2)_\cdot R_\C) &= 
R_\A (u_1, v_2 (B_2 \cap M_{C_2}))_\cdot R_{ \C^{(x_1, x_2)}} \\ 
&=R_\A (u_1, v_2 (B_2 \cap M_{ C_2^{(x_1, x_2)}}))_\cdot 
R_{ \C^{(x_1, x_2)}}\\
&=R_\A(B_1 \times B_2)(u_1, v_2)_\cdot R_{ \C^{(x_1, x_2)}}.
\end{align*}
By \prref{zx}, if 
$R_\A(B_1 \times B_2)(u_1, v_2)_\cdot R_{ \C^{(x_1, x_2)}}
\subset [u_1, v_2]_{\A, \C^{(x_1, x_2)}}/R_{ \C^{(x_1, x_2)}},$
it would follow that 
$R_\A (B_1\times B_2)(v_1, v_2)_\cdot R_\C 
\subset \vvAC/R_\C$.
Consider the sequence of fibrations in \eqref{fibration-sequence}. 
Since
\[
R_\A(B_1 \times B_2)(e, v_2)_\cdot R_{\C^{(\infty)}} 
=R_\A (e, v_2 (B_2 \cap M_{\Cbvv}))_\cdot R_{\C^{(\infty)}} =
\Zcal^{(\infty)},
\]
we see inductively that
$R_\A (B_1\times B_2)(u_1^{(i)}, v_2)_\cdot R_{\C^{(i)}}
\subset \Zcal^{(i)}$
holds for all $i \geq 0$.
\end{proof}

\bpr{pr-shift} Given 
admissible triples $\A = (A_1, A_2, a, K)$
and $\C = (C_1, C_2, c, L)$, 
suppose that $\A^\prime = (A_1, A_2, a, K^\prime)$ and $\C = (C_1, C_2, c, L^\prime)$
are two other admissible quadruples containing the same triples $(A_1, A_2, a)$ and
$(C_1, C_2, c)$. Then there exist $t_2, s_2 \in T_2$
such that
\[
[v_1, v_2]_{\A^\prime, \C^\prime} = (e, t_2) [v_1, v_2]_{\A, \C} (e, s_2),\hs
\forall v_1 \in W_1^{C_1}, v_2 \in {}^{A_2} \! W_2.
\]
\epr

\begin{proof}
Let $K$ be given as in \eqref{K}, and assume that
\[
K^\prime = \{(m_1, m_2) \in M_{A_1} \times M_{A_2} \;| \;\theta_a^\prime(m_1Z_1^\prime) =
m_2Z_2^\prime\},
\]
where, for $i = 1, 2$, $Z_1^\prime$ is a closed subgroup of $Z_{A_i}$, and 
$\theta_1^\prime: M_{A_1}/Z_{1}^{\prime} \to
M_{A_2}/Z_2^\prime$ is an isomorphism having the same properties as $\theta_a$.
The isomorphisms from $M_{A_1}/Z_{A_1}$ to $M_{A_2}/Z_{A_2}$
induced by $\theta_a$ and $\theta_a^\prime$ will still be denoted by the same
symbols. Then our assumptions imply that the automorphism $\theta:=(\theta_a^\prime)^{-1}\theta_a$
of $M_{A_1}/Z_{A_1}$ is inner. Since $\theta$ leaves both $T_1/Z_{A_1}$ and $(B_1 \cap M_{A_1})/Z_{A_1})$
invariant, $\theta = \Ad_{t_1}$ for some $t_1 \in T_1$. It follows that 
$K^\prime(B_1 \times B_2) = (e, t_2)K(B_1 \times B_2)$ for some $t_2 \in T_2$. Similarly, 
$(B_1 \times B_2) L^\prime = (B_1 \times B_2) L (e, s_2)$ for some $s_2 \in T_2$. 
\prref{pr-shift} now follows from  \prref{2dc-plus-2}.
\end{proof}

\begin{example}\label{ex-shift}
Consider the case when $G_1 = G_2 = G$. Take $A_1 = A_2 = \Ga$ and $a = {\rm id}$, where
$\Ga$ is the set of all simple roots for a pair $(B, T)$ of Borel subgroup $B$ of $G$ and
maximal torus $T \subset B$. Take $R_\A = K = G_{{\rm diag}}$ and for some $t \in T$, take
 $R_{\A^\prime}
=K^\prime = \{(g, tgt^{-1}): g \in G\}$. Let  $R_\C =R_{\C^\prime}= B \times B$.
Then it is easy to see that 
\[
R_{\A^\prime} (v, 1) (B \times B) = (e, t)R_\A (v, 1) (B \times B)
\]
for all $v$ in the Weyl group of $(G, T)$.
\end{example}
\sectionnew{Closures of the sets $[v_1, v_2]_{\A, \C}$ in 
$G_1 \times G_2$}\label{closures-GG}
\subsection{The set $\vvAC$ for any $v_1 \in W_1$ and $v_2 \in W_2$}
\label{3begin}
Let $\A = (A_1, A_2, a, K)$ and
$\C = (C_1, C_2, c, L)$ be two arbitrary admissible 
quadruples for $G_1 \times G_2$. Extending 
\eqref{vvBB}, let
\begin{equation}\label{vvAC-general}
\vvAC = R_\A(B_1 \times B_2)(v_1, v_2)(B_1 \times B_2)
R_\C \subset G_1 \times G_2,
\quad \forall
v_1 \in W_1, v_2 \in W_2.
\end{equation}

\ble{lv1}
When $v_1 \in W_1^{C_1}$,
\begin{equation}\label{eq-v1}
(B_1 \times B_2)(v_1, v_2)(B_1 \times B_2)R_\C =
(B_1 \times B_2)(v_1, v_2)R_\C, \hspace{.2in} \forall v_2 \in W_2.
\end{equation}
\ele

\begin{proof}
Since $v_1 \in W_1^{C_1}$, one has  
$B_1 v_1 (B_1 \cap M_{C_1}) = B_1 v_1$, 
and thus
\begin{align*}
(B_1 \times B_2)(v_1, v_2)(B_1 \times B_2)R_\C &=
(B_1 \times B_2)(v_1, v_2)(B_1\cap M_{C_1} \times \{e\})R_\C\\
& =
(B_1 \times B_2)(v_1, v_2)R_\C.X_{(v_1, v_2)} 
\end{align*}
\end{proof}

When 
$v_1 \in W_1^{C_1}$
and $v_2 \in {}^{A_2}\!W_2$, 
the set $\vvAC$ in \eqref{vvAC-general}
is thus the same as what we defined before. Our main result in this section is the
following 
\thref{th-vvAC-closures} which describes 
the closure of $\vvAC$ in $G_1 \times G_2$ for any 
$v_1 \in W_1, v_2 \in W_2$. The proof of \thref{th-vvAC-closures}, which uses
a series of lemmas proved in $\S$\ref{sec-BBRC-closures}, will
be given in $\S$\ref{proof-vvAC}.
In this section, if $X$ is a subset of $G_1 \times G_2$, 
$\overline{X}$
always denotes the closure of $X$ in $G_1 \times G_2$.

\bth{th-vvAC-closures} For any $v_1 \in W_1, v_2 \in W_2$,
\begin{equation}\label{eq-th-vvAC-closures}
\overline{\vvAC} = \bigsqcup_{
\begin{array}{c} v_1^\prime \in W_1^{C_1}, v_2^\prime \in 
{}^{A_2}\!W_2:\\
\exists x_1 \in W_{A_1}, y_1 \in W_{C_1} \; s.t.\\
x_1 v_1^\prime y_1  \leq v_1\\
a(x_1) v_2^\prime c(y_1) \leq v_2\end{array}}
[v_1^\prime, v_2^\prime]_{\A, \C} \hs ({\mbox{disjoint union}}).
\end{equation}
\eth

\bre{re-RABB}
A special case of \eqref{sqcupAC2} says that for any admissible quadruple $\C$, the
$(B_1 \times B_2)$-orbits in $(G_1 \times G_2)/R_\C$ are precisely of the form 
$(B_1 \times B_2) (v_1, v_2)_\cdot R_\C$, where $v_1 \in W_1^{C_1}$ and $v_2 \in W_2$.
Combining \leref{lv1} and \thref{th-vvAC-closures}, we see that 
for any $(B_1 \times B_2)$-orbit $\Ocal$ in 
$(G_1 \times G_2)/R_\C$, the closure
of the set $R_\A \Ocal$ in $(G_1 \times G_2)/R_\C$ 
is a disjoint union of some sets of the form
$R_\A(B_1 \times B_2)(v_1, v_2)_\cdot R_\C$
with
$v_1\in W_1^{C_1}$ and $v_2 \in {}^{A_2}\!W_2$.
\ere
\subsection{The closures of $(B_1 \times B_2, R_\C)$-double cosets in
$G_1 \times G_2$}\label{sec-BBRC-closures} 
In this section, we will describe the 
$(B_1 \times B_2, R_\C)$-double cosets in the set
\[
\overline{(B_1 \times B_2)(v_1, v_2)(B_1 \times B_2)R_\C} 
\subset G_1 \times G_2
\]
for any $v_1 \in W_1, v_2 \in W_2$. Note again that by \leref{lv1}, 
$(B_1 \times B_2)(v_1, v_2)(B_1 \times B_2)R_\C$
is a single $(B_1 \times B_2, R_\C)$-double coset 
when
$v_1 \in W_1^{C_1}$.

For $y_1, z_1 \in W_1$, let
\[
\Wcal_1(y_1, z_1) = \{x_1 \in W_1 \mid B_1x_1B_1 \subset
B_1 y_1 B_1 z_1 B_1\}.
\]
The following \leref{Wyz} on $\Wcal_1(x_1, z_1)$ can be
either proved by induction on the length of $z_1$ or can be seen 
as a direct consequence of the explicit description of 
$\Wcal_1(y_1, z_1)$ in \cite[Remark 3.19]{Borel-Tits}.

\ble{Wyz} Let $y_1, z_1 \in W_1$ be arbitrary. Then every $x_1 
\in \Wcal_1(y_1, z_1)$ is of the form $x_1 = y_1 u_1$ for some
$u_1 \in W_1, u_1 \leq z_1$.
\ele

\ble{BByyRC} Let $y_1 \in W_1, y_2 \in W_2$. Then every
$(B_1 \times B_2, R_\C)$-double coset in 
\[
(B_1 \times B_2)(y_1, y_2)(B_1 \times B_2) R_\C
\]
is of the form $(B_1 \times B_2)(y_1, u_2)R_\C$ 
for some $u_2 \in W_2$,
$u_2 \leq y_2$.
\ele

\begin{proof}
Write $y_2 = w_2 z_2$, where $w_2
\in W_{2}^{C_2}, z_2 \in W_{C_2}$. Then 
\begin{align*}
(B_1 \times B_2)(y_1, y_2)(B_1 \times B_2) R_\C 
&=((B_1 y_1 B_1) \times (B_2 w_2 z_2)) R_\C \\
&=((B_1 y_1 B_1 c^{-1}(z_2^{-1}))\times ( B_2 w_2)) R_\C \\
&=((B_1 y_1 B_1 c^{-1}(z_2^{-1}))\times ( B_2 w_2 B_2)) R_\C \\
&=((B_1 y_1 B_1 c^{-1}(z_2^{-1})B_1) \times (B_2 w_2)) R_\C \\
&=\bigcup_{x_1 \in \Wcal_1(y_1, c^{-1}(z_{2}^{-1}))}((B_1x_1B_1) \times 
(B_2w_2)) R_\C \\
&=\bigcup_{x_1 \in \Wcal_1(y_1, c^{-1}(z_{2}^{-1}))}(B_1 \times 
B_2)(x_1, w_2) R_\C.
\end{align*}
By \leref{Wyz}, every
$x_1 \in \Wcal_1(y_1, c^{-1}(z_2^{-1}))$ is of the form 
$x_1 = y_1 u_1$ for some $u_1 \in W_1$ such that 
$u_1  \leq c^{-1}(z_2^{-1})$,
i.e. $c(u_1^{-1}) \leq z_2$. For such an $x_1 \in  
\Wcal_1(y_1, c^{-1}(z_2^{-1}))$, 
\[
(B_1 \times B_2)(x_1, w_2) R_\C = (B_1\times B_2)
(y_1,  w_2 c(u_1^{-1}))R_\C,
\]
and $w_2 c(u_1^{-1}) \leq w_2 z_2 = y_2$. 
This completes the proof the Lemma.
\end{proof}

\ble{BBRC-closures-0}
For any $v_1 \in W_1$ and $v_2 \in W_2$, one has
\[
\overline{(B_1 \times B_2)(v_1, v_2)(B_1 \times B_2)R_\C} =
\overline{(B_1 \times B_2)(v_1, v_2)(B_1 \times B_2)} \;R_\C.
\]
\ele

\begin{proof} The Lemma follows from \cite[Lemma 2, P. 68]{St} 
(see \leref{lemma-St} in the Appendix) by noting 
that $R_\C/((B_1 \times B_2)\cap R_\C)$ 
is isomorphic to the full flag variety of $M_{C_1}$ (see \leref{Zinfty2}) 
and is hence complete.
\end{proof}

\ble{BBRC-closures-1}
For any $v_1 \in W_1$ and $v_2 \in W_2$, one has
\begin{equation}\label{eq-BBRC-closures-1}
\overline{(B_1 \times B_2)(v_1, v_2)(B_1 \times B_2)R_\C} = \bigcup_{
\begin{array}{l} w_1 \in W_1, w_2\in W_2:\\
w_1 \leq v_1, w_2 \leq v_2\end{array}} (B_1 \times B_2)(w_1, w_2)R_\C.
\end{equation}
\ele

\begin{proof} By \leref{BBRC-closures-0} and the Bruhat decomposition, 
\[
\overline{(B_1 \times B_2)(v_1, v_2)(B_1 \times B_2)R_\C} = 
\bigcup_{\begin{array}{l} y_1 \in W_1, y_2\in W_2:\\
y_1 \leq v_1, y_2 \leq v_2\end{array}} 
(B_1 \times B_2)(y_1, y_2)(B_1 \times B_2)R_\C.
\]
Let $y_1 \in W_1, y_2\in W_2$ be such that 
$y_1 \leq v_1, y_2 \leq v_2$. By \leref{BByyRC}, every 
$(B_1 \times B_2, R_\C)$-double coset in 
$(B_1 \times B_2)(y_1, y_2)(B_1 \times B_2)R_\C$ is of the form 
$(B_1 \times B_2)(y_1, u_2)R_\C$ with 
$u_2 \in W_2, u_2 \leq y_2 \leq v_2$.
Thus $\overline{(B_1 \times B_2)(v_1, v_2)(B_1 \times B_2)R_\C}$ 
is contained in
the right hand side of \eqref{eq-BBRC-closures-1}. Conversely,
let $w_1 \in W_1, w_2\in W_2$ be such that 
$w_1 \leq v_1, w_2 \leq v_2$. Then 
\begin{align*}
(B_1 \times B_2)(w_1, w_2)R_\C &\subset 
(B_1 \times B_2)(w_1, w_2)(B_1 \times B_2)R_\C \\
&\subset
\overline{(B_1 \times B_2)(v_1, v_2)(B_1 \times B_2)}\;R_\C.
\end{align*}
Thus the right hand side of \eqref{eq-BBRC-closures-1} 
is contained in 
$\overline{(B_1 \times B_2)(v_1, v_2)(B_1 \times B_2)R_\C}$.
\end{proof}

\bpr{BBRC-closures1-5}
For any $v_1 \in W_1$ and $v_2 \in W_2$, one has the disjoint union
\[
\overline{(B_1 \times B_2)(v_1, v_2)(B_1 \times B_2)R_\C} = 
\bigsqcup_{\begin{array}{l} w_1 \in W_1^{C_1}, w_2 \in W_2: \\
\exists u_1 \in W_{C_1} \; s.t. \\
w_1 u_1 \leq v_1, w_2 c(u_1) \leq v_2
\end{array}} 
(B_1 \times B_2) (w_1, w_2)R_\C.
\]
\epr

\begin{proof}
\prref{BBRC-closures1-5} 
follows from \leref{BBRC-closures-1} by decomposing
the element $w_1 \in W_1$ in the right hand 
side of \eqref{eq-BBRC-closures-1}
according to the decomposition $W_1 = W_1^{C_1}W_{C_1}$ and by
the fact that 
$(T_1 \times T_2)(u_1, c(u_1))R_\C = (T_1 \times T_2)R_\C$ for 
any $u_1 \in W_{C_1}$, where $T_1$ and $T_2$ are respectively
the maximal tori of $G_1$ and $G_2$ as fixed in \S \ref{quad}. 
The union is disjoint by \thref{LYmain1}.
\end{proof}

\bco{BBRC-closures-2}
For any $v_1 \in W_1^{C_1}, v_2 \in W_2$, one has the disjoint union
\[
\overline{(B_1 \times B_2)(v_1, v_2)R_\C} = \bigsqcup_{
\begin{array}{l} w_1 \in W_1^{C_1}, w_2 \in W_2: \\
\exists u_1 \in W_{C_1} \; s.t. \\
w_1 u_1 \leq v_1, w_2 c(u_1) \leq v_2
\end{array}} 
(B_1 \times B_2) (w_1, w_2)R_\C.
\]
\eco

\bre{get-springer}
By \leref{ww0C} in the Appendix, \coref{BBRC-closures-2} 
is equivalent to
the following description of closures of 
$(B_1 \times B_2, R_\C^-)$-double cosets
in $G_1 \times G_2$ as given in \cite[Lemma 2.2]{Sp}:
for $v_1 \in W_{1}^{C_1},
v_2 \in W_2$,
\[
\overline{(B_1 \times B_2)(v_1, v_2)R_\C^-} = \bigsqcup_{
\begin{array}{l} w_1 \in W_1^{C_1}, w_2 \in W_2: \\
\exists u_1 \in W_{C_1} \; s.t.\\
v_1 u_1^{-1} \leq w_1, w_2 c(u_1) \leq v_2 \end{array}}
(B_1 \times B_2) (w_1, w_2)R_\C^-.
\]
\ere

\ble{wwvv-inter}
Let $w_1 \in W_1, w_2 \in W_2, v_1^\prime \in W_1^{C_1}$, and
$v_2^\prime \in {}^{A_2}\!W_2$. Suppose that
\begin{equation}\label{eq-wwvv-inter}
\left(R_\A (B_1 \times B_2)(w_1, w_2)R_\C\right)
 \cap [v_1^\prime, v_2^\prime]_{\A, \C}
\neq \emptyset.
\end{equation}
Then there exist $x_1 \in W_{A_1}$ and $y_1 \in W_{C_1}$ such that
\[
x_1 v_1^\prime y_1  \leq w_1 \hspace{.2in} {\rm and} \hspace{.2in}
a(x_1) v_2^\prime c(y_1) \leq w_2.
\]
\ele

\begin{proof}
One sees from  \eqref{eq-wwvv-inter} that
$(w_1, w_2)R_\C \subset (B_1 \times B_2)R_\A 
(B_1 \times B_2)(v_1^\prime, v_2^\prime)( B_1 \times B_2)R_\C.$
Since
\[
(B_1 \times B_2)R_\A (B_1 \times B_2) = \bigcup_{x_1 \in W_{A_1}}
(B_1 \times B_2) (x_1, a(x_1))(B_1 \times B_2),
\]
we have
\begin{align*}
(w_1, w_2)R_\C &\subset \bigcup_{x_1 \in W_{A_1}} \left(
(B_1 x_1B_1v_1^\prime B_1) \times 
(B_2 a(x_1) B_2 v_2^\prime B_2)\right) R_\C\\
& = \bigcup_{x_1 \in W_{A_1}} \left(
(B_1 x_1B_1v_1^\prime B_1) \times 
(B_2 a(x_1) v_2^\prime B_2)\right) R_\C,
\end{align*}
where in the last step we used the fact that $l(a(x_1) v_2^\prime)
=l(a(x_1)) + l(v_2^\prime)$.
Thus there exists $x_1 \in W_{A_1}$ such that
$(w_1, w_2) \in ((B_1 x_1B_1v_1^\prime B_1) \times 
(B_2 a(x_1) v_2^\prime B_2)) R_\C.$
Since 
\[
R_\C \subset \bigcup_{y_1 \in W_{C_1}}( (B_1y_1B_1)\times 
(B_2 c(y_1) B_2)),
\]
there exists $y_1 \in W_{C_1}$ such that 
\begin{align*}
(w_1, w_2) &\in (B_1 x_1B_1v_1^\prime B_1y_1B_1) \times 
(B_2 a(x_1) v_2^\prime B_2c(y_1)B_2)\\
&=(B_1 x_1B_1v_1^\prime y_1B_1) \times 
(B_2 a(x_1) v_2^\prime B_2c(y_1)B_2).
\end{align*}
By \leref{Wyz-1} in the Appendix, 
$w_1 \geq x_1 v_1^\prime y_1$ and $ w_2 \geq
a(x_1) v_2^\prime c(y_1)$.
\end{proof}

\subsection{Proof of \thref{th-vvAC-closures}.}\label{proof-vvAC}

Fix $v_1 \in W_1$ and $v_2 \in W_2$. Let 
\[
\Jcal(v_1, v_2) = \{(v_1^\prime, v_2^\prime) 
\in W_1^{C_1} \times {}^{A_2}\!W_2 \mid\;
\overline{\vvAC} \cap 
[v_1^\prime, v_2^\prime]_{\A, \C} \neq \emptyset\}.
\]
Then
\[
\overline{\vvAC} = 
\bigsqcup_{(v_1^\prime, v_2^\prime) \in \Jcal(v_1, v_2)}
\overline{\vvAC} \cap
[v_1^\prime, v_2^\prime]_{\A, \C}.
\]
We will first show that
\begin{align}
\label{calJvv}
\Jcal(v_1, v_2) = \{(v_1^\prime, v_2^\prime) &\in 
W_1^{C_1} \times {}^{A_2}\!W_2
\mid \\
\nonumber
\exists x_1 \in W_{A_1}, y_1 \in W_{C_1} \;  &\mbox{s.t.} \; \;
x_1 v_1^\prime y_1  \leq v_1, \; 
a(x_1) v_2^\prime c(y_1) \leq v_2\}.
\end{align}
Indeed, by \leref{lemma-St} in the Appendix,  
\[
\overline{\vvAC} = R_\A \;
\overline{(B_1 \times B_2)(v_1, v_2)(B_1 \times B_2) R_\C}.
\]
Thus by \leref{BBRC-closures-1}, 
\[
\overline{\vvAC} = \bigsqcup_{
\begin{array}{l} w_1 \in W_1, w_2\in W_2:\\
w_1 \leq v_1, w_2 \leq v_2\end{array}} 
R_\A(B_1 \times B_2)(w_1, w_2)R_\C.
\]
Suppose that
$(v_1^\prime, v_2^\prime) \in \Jcal(v_1, v_2)$. Then 
there exists $(w_1, w_2) \in W_1 \times W_2$ with 
$w_1 \leq v_1, w_2 \leq v_2$
such that
\[
(R_\A(B_1 \times B_2)(w_1, w_2)R_\C) \cap  
[v_1^\prime, v_2^\prime]_{\A, \C}
\neq \emptyset.
\]
By \leref{wwvv-inter}, there exist 
$x_1 \in W_{A_1}, y_1 \in W_{C_1}$
such that 
\[
x_1 v_1^\prime y_1  \leq w_1 \leq v_1 \hspace{.2in} 
{\rm and} \hspace{.2in} 
a(x_1) v_2^\prime c(y_1) \leq w_2 \leq v_2.
\]
Thus $(v_1^\prime, v_2^\prime)$
is in the set of the right hand side of \eqref{calJvv}. 
Conversely, suppose 
that $(v_1^\prime, v_2^\prime) \in W_1^{C_1} \times {}^{A_2}\!W_2$
are such that $x_1 v_1^\prime y_1 \leq v_1$ and 
$a(x_1) v_2^\prime c(y_1) \leq v_2$ for some $x_1 \in W_{A_1}$
and $y_1 \in W_{C_1}$. Let $w_1 = x_1 v_1^\prime y_1$ and 
$w_2 = a(x_1) v_2^\prime c(y_1)$ so that
\[
v_1^\prime = x_1^{-1} w_1 y_1^{-1} \hspace{.2in}
{\rm and} \hspace{.2in} 
v_2^\prime = a(x_1^{-1}) w_2 c(y_1^{-1}).
\]
It follows that
$(v_1^\prime T_1, \; v_2^\prime T_2) = (x_1^{-1} w_1 y_1^{-1} T_1, \; 
a(x_1^{-1}) w_2 c(y_1^{-1})T_2)$ and hence
\begin{equation}\label{vTvT}
(v_1^\prime T_1, \; v_2^\prime T_2)
 \subset R_\A (B_1 \times B_2)(w_1, w_2)R_\C
\subset \overline{\vvAC},
\end{equation}
where $T_i = B_i^- \cap B_i$ for $i = 1, 2$. 
Thus $\overline{\vvAC} \cap [v_1^\prime, \; v_2^\prime]_{\A, \C}
\neq \emptyset$, and so $(v_1^\prime, \; v_2^\prime) \in \Jcal(v_1, v_2)$.
\qed

For any $(v_1^\prime, \; v_2^\prime) \in \Jcal(v_1, v_2)$,
it  follows from \eqref{vTvT} that 
\[
R_\A (v_1^\prime, \; v_2^\prime T_2) R_\C
=R_\A (v_1^\prime T_1, \; v_2^\prime T_2) R_\C \subset \overline{\vvAC}.
\]
By \leref{RvvT}, $R_\A (v_1^\prime, \; v_2^\prime T_2) R_\C$ is 
dense in $[v_1^\prime, \; v_2^\prime]_{\A, \C}$.
Hence $[v_1^\prime, \; v_2^\prime]_{\A, \C} \subset \overline{\vvAC}$,
and $\overline{\vvAC} \cap [v_1^\prime, v_2^\prime]_{\A, \C} 
=[v_1^\prime, v_2^\prime]_{\A, \C}$.
This completes the proof of \thref{th-vvAC-closures}.
\sectionnew{The $(R_\A, R_\C^-)$-stable subsets $\vvACm$}\label{minus-strata}

\subsection{The subsets $\vvACm$ of $G_1 \times G_2$}
\label{vvACm-def}
Let again $\A$ and $\C$ be two admissible quadruples for $G_1 \times G_2$.
Let $w_{0,\Ga_1}$ and $\wc$ be the longest element in $W_1$ and $W_{C_1}$
respectively. Associated to $\C$, we have the admissible quadruple
$\C^* \stackrel{{\rm def}}{=} (C_{1}^{*}, C_{2}^{*}, c^{*}, L^*),$
where
\begin{equation}\label{Cstar}
C_{1}^{*} = -\wg(C_1), \hspace{.1in} C_2^* = C_2, \hspace{.1in}
c^* = c(\wg\wc)^{-1}, \hspace{.1in} 
L^* = \Ad_{(\dwg \dwc, \, e)}L.
\end{equation}
It is easy to see that $R_{\C^*} = \Ad_{(\dwg \dwc, \, e)} R_\C^-.$
Moreover, 
\begin{equation}\label{WC}
W_{1}^{C_1} (\wg\wc)^{-1} = W_{1}^{C_1^*}
\end{equation}

\bpr{vvACm-strata}
For $v_1 \in W_{1}^{C_1}$ and $v_2 \in {}^{A_2}\! W_2$, let
\[
\vvACm = R_\A(v_1, v_2M_{\Cbvv})R_\C^-
\subset G_1 \times G_2.
\]
Then the following properties hold:

(i)
\[
G_1 \times G_2= \bigsqcup_{v_1 \in W_{1}^{C_1}, 
v_2 \in {}^{A_2}\! W_2} \vvACm \hspace{.3in} (\mbox{disjoint union});
\]
 
(ii) $\vvACm = R_\A(B_1 \times B_2)(v_1, v_2) R_\C^-$ for
every $v_1 \in W_{1}^{C_1}$ and $v_2 \in {}^{A_2}\! W_2$;

(iii) $[v_1, v_2]_{\A, \C}$ is locally closed, smooth, and 
irreducible.
Its projection $[v_1, v_2]_{\A, \C}/R_\C^-$ to 
$(G_1 \times G_2) /R_\C^-$ fibers over the
flag variety $M_{A_1}/(M_{A_1} \cap P_{\Aavv})$
with fibers isomorphic to the product of $M_{\Cbvv}/Y_2$ 
and the affine space of 
dimension $\dim U_{A_1(v_1, v_2)} - l(v_1) + l(v_2)$,
where $Y_2$ is as in \thref{main1}.
\epr

\begin{proof} 
Let $v_1^* = v_1(\wg \wc)^{-1} \in W_1^{C_1^*}$. All the statements in 
\prref{vvACm-strata} follow from the fact that 
\[
\vvACm = R_\A(B_1 \times B_2)(v_1^*, v_2)
(B_1 \times B_2)R_{\C^*}(\dwg \dwc, e).
\]
\end{proof}

\subsection{Closures of the sets $\vvACm$ in
$G_1 \times G_2$}\label{closures-RARCminus}

For each $v_1 \in W_1$ and $v_2 \in W_2$,
set
\begin{equation}\label{eq-vvACm}
\vvACm = R_\A(B_1 \times B_2)(v_1, v_2)(w_{0, C_1}(B_1^-) \times B_2)
R_\C^- \subset G_1 \times G_2.
\end{equation}
It follows from \leref{lv1} that  
\begin{equation}\label{bb-minus}
\vvACm = R_\A(B_1 \times B_2)(v_1, v_2)R_\C^-, \hspace{.2in}
\mbox{when} \; \; v_1 \in W_{1}^{C_1}, \, v_2 \in W_2.
\end{equation}

\bth{th-vvACm-closures}
For any $v_1 \in W_1$ and $v_2 \in W_2$, with $\vvACm$ given in \eqref{eq-vvACm},
one has
\[
\overline{\vvACm} = \bigsqcup_{\begin{array}{c} 
v_1^\prime \in W_{1}^{C_1}, v_2^\prime \in {}^{A_2}\! W_2:\\
\exists x_1 \in W_{A_1}, y_1 \in W_{C_1} \; s.t.\\
x_1 v_1^{\prime} y_1 \wc \geq v_1\wc  \\
a(x_1) v_2^\prime c(y_1) \leq v_2\end{array}}
[v_1^\prime, v_2^\prime]_{\A, \C}^-.
\]
\eth

\begin{proof}
For $v_1 \in W_1$, let again 
$v_1^* = v_1(\wg \wc)^{-1}$. Then
\[
\vvACm = R_\A(B_1 \times B_2)(v_1^*, v_2)
(B_1 \times B_2)R_{\C^*}(\dwg \dwc, e),
\hspace{.2in} \forall v_1 \in W_1, \; v_2 \in W_2.
\]
For
$y_1 \in W_{C_1}$, set $(y_1)_* = 
(\wg w_{0, C_1}) y_1 (\wg w_{0, C_1})^{-1} \in W_{C_1^*}$.
By \thref{th-vvAC-closures}, 
\begin{align*}
\overline{\vvACm} &= \bigsqcup_{\begin{array}{c} 
v_1^\prime \in W_{1}^{C_1}, v_2^\prime \in {}^{A_2}\! W_2:\\
\exists x_1 \in W_{A_1}, y_1 \in W_{C_1} \; s.t.\\
x_1 (v_1^\prime)^* (y_1)_* \leq v_1^*\\
a(x_1) v_2^\prime c(y_1) \leq v_2\end{array}}
[v_1^\prime, v_2^\prime]_{\A, \C}^-\\
&=\bigsqcup_{\begin{array}{c} 
v_1^\prime \in W_{1}^{C_1}, v_2^\prime \in {}^{A_2}\! W_2:\\
\exists x_1 \in W_{A_1}, y_1 \in W_{C_1} \; s.t.\\
x_1 v_1^{\prime} y_1 \wc \geq v_1\wc  \\
a(x_1) v_2^\prime c(y_1) \leq v_2\end{array}}
[v_1^\prime, v_2^\prime]_{\A, \C}^-.
\end{align*}
\end{proof}
\sectionnew{The $R_\A$-stable pieces in $\ol{G}$ and their closures}
\label{proof_th}
We retain the notation in \S \ref{wond}. In particular, 
$G$ is a semi-simple algebraic group of adjoint type, and 
$\overline{G}$ denotes the De Concini-Procesi compactification of $G$.
In this section, unless otherwise stated, if $X$ is a subset of 
$\overline{G}$,
$\overline{X}$ always denotes the closure of $X$ in $\overline{G}$.

For each $J \subset G$, recall that
$h_J$ is a point in $\overline{G}$ such that the
stabilizer subgroup of $G \times G$ at $h_J$ is $R_J^-$
given in \eqref{RJ-minus}. Let $\A =
(A_1, A_2, a, K)$ be any admissible quadruple for $G \times G$. Then 
we have the decomposition of $\ol{G}$ into $R_\A$-stable pieces
\[
\ol{G} = \bigsqcup_{J \subset \Ga, v_1 \in W^J,
v_2 \in {}^{A_2} \! W} [J, v_1, v_2]_\A,
\]
where for $J \subset \Ga$ and for $v_1 \in W^J$ and $v_2 \in {}^{A_2}\!W$,
the subset $[J, v_1, v_2]_\A$ of $\overline{G}$ 
is defined by \eqref{wond_part},
namely 
\[
[J, v_1, v_2]_\A = R_\A(B \times B)(v_1, v_2)_\cdot h_J.
\]
By part (iii) of \prref{vvACm-strata}, we have 

\bpr{pr-JvvA}
Let $A$ be any admissible quadruple for $G \times G$. Then for any $J \subset \Ga$
and $v_1 \in W^J, v_2 \in {}^{A_2}\!W$, $[J, v_1, v_2]_\A$ is
a locally closed smooth subset of $\ol{G}$. It fibers over the
flag variety $M_{A_1}/(M_{A_1} \cap P_{\Aavv})$
with fibers isomorphic to the product of $M_{J(v_1, v_2)}/Z_J$ 
and the affine space of 
dimension $\dim U_{A_1(v_1, v_2)} - l(v_1) + l(v_2)$, where 
$J(v_1, v_2)$ is the smallest subset of $J$ stable under
$v_2^{-1}av_1$ and $A_1(v_1, v_2) = v_1J(v_1, v_2)\subset A_1$.
\epr

When $R_\A = G_{{\rm diag}}$, the \prref{pr-JvvA} coincides with 
the description of the geometry of Lusztig's $G_{{\rm diag}}$-stable
pieces given in \cite{H2}. 

In this section, we study the closures of the subsets $[J, v_1, v_2]_\A$ in $\ol{G}$.

\subsection{The first description}\label{sec-first-descrip}
For $J \subset \Ga$ and $v_1 \in W^J$, $v_2 \in W$, let
\[
[J, v_1, v_2] = (B \times B)(v_1, v_2)_\cdot h_J \subset \ol{G}.
\] 
The following \leref{clos1} follows immediately from \leref{lemma-St} in 
the Appendix.

\ble{clos1}
For any admissible triple $\A$ for 
$G \times G$ and all $J \subset \Gamma,$
$v_1 \in W^J,$ $v_2 \in {}^{A_2} \! W$,
\[
\ol{[J, v_1, v_2]_\A} = R_\A \ol{[J, v_1, v_2]}.
\]
\ele

\bpr{withw0C} Let $\A$ be 
an admissible quadruple for $G \times G$. Then
for any $J \subset \Ga$ and $(v_1, v_2) \in W^J \times {}^{A_2}\!W$,
\begin{equation}\label{eq-withw0C}
\ol{[J, v_1, v_2]_\A} = \bigsqcup_{I \subset J}
\bigsqcup_{\begin{array}{c}v_1^\prime \in W^{I}, v_2^\prime 
\in {}^{A_2}\! W:\\
\exists x \in W_{A_1}, y \in W_I, z \in W_J^I \; s.t.\\
l(v_2z) = l(v_2) + l(z)\\
x v_1^{\prime} y w_{0, I} \geq v_1z w_{0, I}  \\
a(x) v_2^\prime y \leq v_2z
\end{array}} [I, \; v_1^\prime, \; v_2^\prime]_{\A}.
\end{equation}
\epr

\begin{proof}
It is well-known \cite{DP} that
for $I_1, I_2 \subset \Ga$, 
\[
\overline{(G \times G)_\cdot h_{I_1}} \; \cap\; 
(G \times G)_\cdot h_{I_2} \neq 
\emptyset
\]
if and only if $I_2 \subset I_1$. Thus 
\[
\ol{[J, v_1, v_2]}_\A = \bigsqcup_{I \subset J} 
\left(\ol{[J, v_1, v_2]_\A}
\;\cap\; (G \times G)_\cdot h_I\right).
\]
Fix $I \subset J$. By \leref{clos1}, 
\[
\ol{[J, v_1, v_2]_\A}
\;\cap \;(G \times G)_\cdot h_I = R_\A \left(\ol{[J, v_1, v_2]}
\;\cap \;(G \times G)_\cdot h_I\right).
\]
By \cite[Lemma 2.3]{Sp},
\[
\ol{[J, v_1, v_2]}
\;\cap \;(G \times G)_\cdot h_I = 
\bigsqcup_{\begin{array}{c}z \in W_{J}^{I},\\
l(v_2z) = 
l(v_2) + l(z)\end{array}} \ol{[I,\;  v_1z, \; v_2z]}^{I},
\]
where, for a subsets $Y$ of $(G \times G)_\cdot h_I$, $\ol{Y}^I$ 
denotes the
closure of $Y$ in $(G \times G)_\cdot h_I$. Thus 
\[
\ol{[J, v_1, v_2]_\A}
\;\cap \; (G \times G)_\cdot h_I = 
\bigsqcup_{\begin{array}{c}z \in W_{J}^{I}, \\
l(v_2z) = 
l(v_2) + l(z)\end{array}} R_\A \; \ol{[I, \; v_1z, \; v_2z]}^{I}.
\]
By \leref{lemma-St}, $R_\A \; \ol{[I, \; v_1z, \; v_2z]}^{I} =
\ol{R_\A[I, \; v_1z, \; v_2z]}^{I}$ for every $z \in W_J^I$. The
decomposition in \eqref{eq-withw0C} 
now follows from \thref{th-vvACm-closures}.
\end{proof}

In the following $\S$\ref{sec-second-descrip} and $\S$\ref{sec-third-descrip}, we will 
simplify the descriptions of the the closure relations in \prref{withw0C}.

\subsection{The second description}\label{sec-second-descrip}

Recall that $(G \times G)_\cdot h_\emptyset \cong (G \times G)/(B_- \times B)$
is the unique closed $(G \times G)$-orbit in $\ol{G}$. For 
$J \subset \Ga$ and $v_1 \in W^J, \, v_2  \in{}^{A_2} \! W$,
let 
\[
\partial_{\emptyset}\overline{[J, v_1, v_2]_\A}
= \overline{[J, v_1, v_2]_\A} \cap (G \times G)_\cdot h_{\emptyset}.
\]

\bth{th-boundary}
Let $\A$ be an admissible quadruple for 
$G \times G$. Let $I, J \subset \Gamma,$
$v_1 \in W^J, v_1^\prime \in W^I$, and 
$v_2, v_2^\prime \in  {}^{A_2} \! W$. Then the following are equivalent:

(i) $[I, v_1^\prime, v_2^\prime]_\A \subset \overline{[J, v_1, v_2]_\A}$;

(ii) $I \subset J$ and $\partial_\emptyset \overline{[I, v_1^\prime, v_2^\prime]_\A} 
\subset \partial_\emptyset \overline{[J, v_1, v_2]_\A}$;

(iii) $I \subset J$ and $[\emptyset, v_1^\prime, v_2^\prime]_\A \subset
\partial_\emptyset \overline{[J, v_1, v_2]_\A}$.
\eth

\begin{proof}
Clearly (i) implies (ii). Since $[\emptyset, v_1^\prime, v_2^\prime]_{\A}
\subset \partial_\emptyset \overline{[I, v_1^\prime, v_2^\prime]_\A}$,  
(ii) implies (iii). It remains to show that (iii) implies (i).

Assume (iii). By \prref{withw0C}, there exist $x \in W_{A_1}$ and $z \in W_J$ with
$l(v_2 z) = l(v_2) + l(z)$ such that
$xv_1^\prime \geq v_1z$ and $a(x)v_2^\prime \leq v_2z$. Write
$z = u y$ with $u \in W^I_J$ and $y \in W_I$. Then the set
\[
S = \{(x^\prime, y^\prime) \in W_{A_1} \times W_I \mid 
x^\prime v_1^\prime \geq v_1uy^\prime, \; \;
a(x^\prime)v_2^\prime \leq v_2uy^\prime\}
\]
is non-empty. Let $(x_0, y_0) \in W_{A_1} \times W_I$ be a 
minimal element in $S$. We claim that 
\begin{equation}\label{x0y0}
(x_0v_1^\prime)^I \geq v_1 u y_0 (x_0 v_1^\prime)_I^{-1} 
\hspace{.2in} 
{\rm and} \hspace{.2in} a(x_0) v_2^\prime y_0^{-1} \leq v_2 u,
\end{equation}
where $(x_0v_1^\prime)^I \in W^I$ and $(x_0v_1^\prime)_I \in W_I$ are
such that $x_0 v_1^\prime = (x_0v_1^\prime)^I (x_0v_1^\prime)_I$.
By \leref{ww0C}, it would follow from \eqref{x0y0} that
\[
x_0 v_1^\prime y_0^{-1} w_{0, I} \geq v_1 u w_{0,I}
\hspace{.2in} {\rm and} \hspace{.2in} 
a(x_0) v_2^\prime y_0^{-1} \leq v_2u,
\]
and, since $l(v_2u) = l(v_2) + l(u)$,  
we would see by \prref{withw0C} that  
$[I, v_1^\prime, v_2^\prime]_\A
\subset \overline{[J, v_1, v_2]_\A}$.

It remains to prove \eqref{x0y0}. We first show that 
$a(x_0) v_2^\prime y_0^{-1}\leq v_2u$. 
Indeed, since $a(x_0)v_2^\prime \leq  v_2uy_0$, it follows from
\leref{wwu} in the Appendix that there exists $y_1 \leq y_0$ such that
\begin{equation}\label{v2uy1}
v_2u = v_2 u y_0 y_0^{-1} \geq a(x_0) v_2^\prime y_1^{-1}.
\end{equation}
Again by \leref{wwu}, there exists $y_2 \leq y_1$ such that 
\[
a(x_0) v_2^\prime = a(x_0) v_2^\prime y_1^{-1} y_1 \leq v_2u y_2.
\]
Since $y_2 \leq y_1 \leq y_0$, we have 
$x_0 v_1^\prime \geq v_1 u y_0 \geq v_1 u y_2$.
Thus $(x_0, y_2) \in S$. Since $(x_0, y_2) \leq (x_0, y_0)$ and since
$(x_0, y_0)$ is minimal in $S$, we must have $y_2 = y_0$. Hence $y_1 = 
y_0$, and $a(x_0) v_2^\prime y_0^{-1}\leq v_2u$ by \eqref{v2uy1}.

We now show that $l(x_0 v_1^\prime) = l(x_0) + l(v_1^\prime)$. If
$l(x_0 v_1^\prime) < l(x_0) + l(v_1^\prime)$, then by \leref{uwvw} 
in the Appendix, there exists $x_1 < x_0$ such that 
$x_1 v_1^\prime > x_0 v_1^\prime \geq v_1 u y_0$. Since $a(x_1)v_2^\prime 
\leq a(x_0) v_2^\prime \leq v_2 u y_0$, we have $(x_1, y_0) \in S$. Since
$(x_1, y_0) < (x_0, y_0)$, this is a contradiction to the minimality of
$(x_0, y_0)$ in $S$. Hence $l(x_0 v_1^\prime) = l(x_0) + l(v_1^\prime)$.

By \leref{wwu} in the Appendix, there exists $y_1 \leq (x_1v_1^\prime)_I$
such that
\[
(x_0v_1^\prime)^I = (x_0v_1^\prime) (x_0v_1^\prime)_I^{-1} 
\geq v_1 u y_0 y_1^{-1}.
\]
By \leref{uwxv} in the Appendix, there exists $x_2 \leq x_0$ such that
$(x_0v_1^\prime)^I y_1 = x_2 v_1^\prime$.  Now 
\[
x_2 v_1^\prime = (x_0v_1^\prime)^I y_1 \geq 
v_1 u y_0 y_1^{-1} y_1 = v_1uy_0,
\]
and $a(x_2)v_2^\prime \leq a(x_0) v_2^\prime \leq v_2 u y_0$. Hence 
$(x_2, y_0) \in S$. By the minimality of $(x_0, y_0)$ in $S$, we have
$x_2 = x_0$, so $y_1 = (x_1v_1^\prime)_I$, and 
$(x_0v_1^\prime)^I \geq v_1 u y_0 (x_1v_1^\prime)_I^{-1}$. This proves
\eqref{x0y0}.
\end{proof}

As a corollary of (iii) in \thref{th-boundary} and \prref{withw0C}, we get the following
second description of the closure relations on the $R_\A$-stable pieces in $\ol{G}$.

\bco{co-second-descrip} 
Let $\A$ be an admissible quadruple for $G \times G$. Then
for any $J \subset \Ga$ and $(v_1, v_2) \in W^J \times {}^{A_2}\!W$,
\begin{equation}\label{eq-second}
\ol{[J, v_1, v_2]_\A} = \bigsqcup_{I \subset J}
\bigsqcup_{\begin{array}{c}v_1^\prime \in W^{I}, v_2^\prime 
\in {}^{A_2}\! W:\\
\exists x \in W_{A_1}, z \in W_J \; s.t.\\
l(v_2z) = l(v_2) + l(z)\\
x v_1^{\prime} \geq v_1z, \; \; 
a(x) v_2^\prime \leq v_2z
\end{array}} [I, \; v_1^\prime, \; v_2^\prime]_{\A}.
\end{equation}
\eco

\subsection{The third description}\label{sec-third-descrip} In this section, we show that
the condition $l(v_2z) = l(v_2) + l(z)$ in \eqref{eq-second} can be dropped. 

\bth{th-third-descrip} Let $\A$ be an admissible quadruple for $G \times G$. Then
for any $J \subset \Ga$ and $(v_1, v_2) \in W^J \times {}^{A_2}\!W$,
\begin{equation}\label{eq-third}
\ol{[J, v_1, v_2]_\A} = \bigsqcup_{I \subset J}
\bigsqcup_{\begin{array}{c}v_1^\prime \in W^{I}, v_2^\prime 
\in {}^{A_2}\! W:\\
\exists x \in W_{A_1}, z \in W_J \; s.t.\\
x v_1^{\prime} \geq v_1z, \; \; 
a(x) v_2^\prime \leq v_2z
\end{array}} [I, \; v_1^\prime, \; v_2^\prime]_{\A}.
\end{equation}
\eth

\begin{proof} Fix $J \subset \Ga$ and $(v_1, v_2) \in W^J \times {}^{A_2}\!W$.
It is enough to show that the right hand side of \eqref{eq-third} is contained in
the right hand side of \eqref{eq-second}. To this end, let
$I \subset J$ and let $(v_1^\prime, v_2^\prime)  \in W^{I} \times
{}^{A_2}\! W$ be such that there exist $x \in W_{A_1}$ and $z \in W_J$ with
$x v_1^{\prime} \geq v_1z$ and $a(x) v_2^\prime \leq v_2z$. Choose such an $x \in W_{A_1}$ and let
\[
Z = \{z^\prime \in W_J\; | \; \; x v_1^{\prime} \geq v_1z^\prime, \; \; a(x) v_2^\prime \leq v_2z^\prime\}.
\]
Then $Z \neq \emptyset$. Let $z_0 \in Z$ be a minimal element. We claim that 
$l(v_2z_0) = l(v_2) + l(z_0)$. Indeed, if $l(v_2z_0) < l(v_2) + l(z_0)$, then by \leref{uwvw}
in the Appendix, there exists $z_1 < z_0$ such that $v_2 z_1 > v_2 z_0$. Thus 
$v_2 z_1 > a(x) v_2^\prime$. Since $v_1 \in W^J$ and $z_1, z_0 \in W_J$, we also have 
$v_1 z_1 < v_1 z_0 \leq x v_1^\prime$. Thus $z_1 \in Z$, which contradicts to the
fact the $z_0$ is a minimal element in $Z$. This shows that  $l(v_2z_0) = l(v_2) + l(z_0)$,
and thus $[I, v_1^\prime, v_2^\prime]_\A$ is contained in the right hand side of \eqref{eq-second}.
\end{proof}

\bre{re-v1} 
Note that in the proofs of \prref{withw0C}, \coref{co-second-descrip}, and
\thref{th-third-descrip}, 
we did not use the fact that
$v_2 \in {}^{A_2}\!W$. In fact, the decomposition formulas \eqref{eq-withw0C}, \eqref{eq-second},
and \eqref{eq-third}
hold for any $v_1 \in W^J$ and $v_2 \in W_2$. Thus 
the closure of $R_\A \Ocal$ in $\overline{G}$ for any $(B \times B)$-orbit $\Ocal$ in 
$\ol{G}$ is a union of the sets of the form 
$[I, v_1^\prime, v_2^\prime]_\A$ for $I \subset \Ga$,
$v_1^\prime \in W^I$, and $v_2^\prime \in {}^{A_2}\!W$. \prref{withw0C}, \coref{co-second-descrip}, and
\thref{th-third-descrip} give three equivalent
descriptions of the decomposition. See also \reref{re-RABB}.
\ere

\sectionnew{Appendix}\label{appendices}

\subsection{A lemma from \cite{St}}\label{appendix-lemma-St}

The following is \cite[Lemma 2, P. 68]{St}. We state it here for the 
convenience of the reader.

\ble{lemma-St}
Let $G$ be an algebraic group acting on a variety $V$. Let
$H$ be a closed subgroup of $G$ and let $U \subset V$ be a 
closet subset of $V$, invariant under the action of $H$.
Assume that $G/H$ is complete. Then $G_\cdot U$ is closed.
\ele

\subsection{A few facts on Weyl groups}\label{bruhat}
Let $G$ be a any connected reductive algebraic group over
an algebraically closed field. Let $T$ be a maximal torus
of $G$, let $B$ a Borel subgroup of $G$ containing $T$, and let
$\Ga$ be the set of simple roots
for $(B, T)$.
Let $W$ be the
Weyl group of $\Ga$. For the convenience of the reader, we collect 
in this section a few facts on $W$ that
are used in this paper. Recall that for two 
subsets $A$ and $C$ of $\Gamma$ we denote 
${}^{A}\!W^C = {}^{A}\!W \cap W^C$.

\ble{car-prop275}\cite[Proposition 2.7.5]{Car}
For any $A, C \subset \Gamma$, every $v \in {}^{A}\!W$ can be uniquely
written as a product $v = xu$, where $x \in {}^{A}\!W^C$ and
$u \in {}^{\!C \cap x^{-1}(A)}W_{C}$.
\ele

For $C \subset \Ga$, let $w_{0, C}$ be the longest element of $W_C$.

\ble{ww0C}\cite[Page 79]{Sp}
Let $x, w \in W^C$ and $u \in W_C$. Then $xu^{-1} \leq w$  if and only
if $x \wcc \leq w u \wcc$.
\ele

\begin{proof}
If $x u^{-1} \leq w$, then $x \wcc = xu^{-1} u \wcc \leq w u \wcc$.
Conversely, assume that $x \wcc \leq w u \wcc$. 
Then there exist $w_1 \leq w$ and
$y \leq u \wcc$ such that $x \wcc = w_1 y$. 
It follows from $y \leq u \wcc$ that
$y \wcc \geq u$, so $u^{-1} \leq \wcc y^{-1}$.  Thus
$x u^{-1} \leq x \wcc y^{-1} = w_1 \leq w.$
\end{proof} 

For $y, z \in W$, let
$\Wcal(y, z) =\{x \in  W \mid BxB \subset ByBzB\}.$

\ble{Wyz-1} \cite[p. 420]{Sp-lect}
Let $y, z \in W$. Then $x \geq yz$ 
for any $x \in \Wcal(y, z)$.
\ele

\ble{uwmax} \cite[Lemma 3.3]{H2}
For any $u, w \in W$, the subset $\{vw \mid v \leq u\}$ 
of $W$ contains a 
unique maximal element $u_1w$. Moreover, $l(u_1 w) = l(u_1) + l(w)$.
\ele

\ble{uwvw}
If $u, w \in W$ are such that $l(uw) < l(u) + l(w)$, then there exists
$u_1$ such that $u_1 < u$ and $u_1 w > uw$.
\ele

\begin{proof} The $u_1$ such that $u_1w$ is the maximal element
in the set $\{vw \mid v \leq u\}$ is as required.
\end{proof}

\ble{wwu} \cite[Corollary 3.4]{H2}
Let $u, w, w^\prime \in W$ and assume that $w^\prime \leq w$. Then

(i) there exists $u_1 \leq u$ such that $w^\prime u_1 \leq wu$;

(ii) there exists $u_2 \leq u$ such that $w^\prime u \leq wu_2$.
\ele

\ble{uwxv}\cite[Lemma 3.10]{H2}
Let $J \subset \Ga, w \in W^J$, and $u \in W$ be such that 
$l(uw) = l(u) + l(w)$.
Write $uw = xv$ with $x \in W^J$ and $v \in W_J$. 
Then for any $v^\prime \leq v$,
there exists $u^\prime \leq u$ such that $u^\prime w = x v^\prime$.
\ele

\end{document}